\documentclass[preprint,12pt]{elsarticle}

\usepackage[top=1in,bottom=1in,left=1in,right=1in]{geometry}

\usepackage[utf8]{inputenc}
\usepackage[T1]{fontenc}

\usepackage{amssymb}
\usepackage{amsmath}
\usepackage{mathtools}
\usepackage{bm}

\usepackage{graphicx}
\usepackage{float}
\usepackage{booktabs}
\usepackage{adjustbox}
\usepackage{multirow}

\usepackage{algpseudocode}
\usepackage[ruled,vlined]{algorithm2e}

\usepackage{enumitem}
\usepackage{nicefrac}
\usepackage{url}
\usepackage{microtype}

\usepackage{xcolor}
\usepackage{hyperref}
\hypersetup{colorlinks=true,linkcolor=blue,citecolor=blue,urlcolor=blue}

\usepackage{amsmath,amsfonts,bm,mathtools,amsthm,amssymb}

\newcommand{\bsmat}{\left[\begin{smallmatrix}}
\newcommand{\esmat}{\end{smallmatrix}\right]}

\newcommand{\R}{\mathbb{R}}

\DeclareMathAlphabet{\mathsfit}{\encodingdefault}{\sfdefault}{m}{sl}
\SetMathAlphabet{\mathsfit}{bold}{\encodingdefault}{\sfdefault}{bx}{n}

\newcommand{\vare}{\varepsilon}

\newtheorem{theorem}{Theorem}[section]
\newtheorem{proposition}{Proposition}[section]
\newtheorem{definition}{Definition}[section]

\newtheorem{lemma}{Lemma}[section]

\newtheorem{corollary}{Corollary}[section]

\biboptions{sort&compress}

\myfooterfont{\normalfont\footnotesize}
\myfooter[C]{\parbox{0.98\textwidth}{\centering
  This paper has been accepted at NeurIPS 2026.\\[3pt]
  \thepage}}

\begin{document}

\begin{frontmatter}

\title{Physics-Informed Optimal Control by Control-Only Supervision with Error Guarantees on Value and Policy}

\author[inst1,inst2,inst3]{Zihua Wang}
\author[inst4]{Xiaopei Jiao\corref{cor1}}
\author[inst1]{Yunfeng Cai}

\cortext[cor1]{\raggedright Corresponding author. \\
E-mails:
\texttt{wangzihua@bimsa.cn}, \texttt{jiaoxiaopei@cnu.edu.cn},
\texttt{caiyunfeng@bimsa.cn}.}

\address[inst1]{Beijing Institute of Mathematical Sciences and Applications (BIMSA), Beijing, China}
\address[inst2]{State Key Laboratory of Mathematical Science, Academy of Mathematics and Systems Science, Chinese Academy of Sciences, Beijing, China}
\address[inst3]{University of Chinese Academy of Sciences, Beijing, China}
\address[inst4]{Capital Normal University, Beijing, China}

\begin{abstract}
We propose a physics-informed method for solving high-dimensional optimal control problems that requires only control labels, not value function labels. A neural network parameterizes the value function, and the control policy is derived from its gradient. The network is trained by minimizing a weighted sum of the standard Hamilton-Jacobi-Bellman (HJB) residual and the mismatch between the network’s control output and control labels, which are precomputed by a state-dependent Riccati equation (SDRE) solver. Theoretically, with the help of the viscous HJB approach, we establish error bounds for the value function and its gradient by applying stochastic dynamic programming and Liouville-type inequality, respectively, providing conditional error guarantees under the stated assumptions.
Numerical experiments show that our neural surrogate model not only captures the true control policy exactly but also closely matches the cost-to-go of the SDRE reference with significant inference speedup and excellent stability. Ablation studies fully exhibit the merits of our algorithm in noise robustness, out-of-distribution (OOD) generalization, etc.
\end{abstract}

\begin{keyword}
physics-informed neural networks \sep optimal control \sep
Hamilton--Jacobi--Bellman equation \sep state-dependent Riccati equation \sep
error analysis \sep viscosity solution
\end{keyword}

\end{frontmatter}

\section{Introduction}
Optimal control problem (OCP) provides a systematic framework for determining control laws that minimize a specified performance criterion while satisfying system dynamics and constraints. By formulating control objectives as mathematical optimization problems, this approach enables the design of controllers that achieve desired behaviors such as trajectory tracking \cite{peng2019symplectic}, disturbance rejection \cite{vidyasagar2003optimal}, or minimal energy consumption \cite{sardellitti2011optimal} in an optimal sense. The resulting control policies balance competing objectives, such as speed of response and control effort, and are fundamental to applications ranging from aerospace and robotics to economics \cite{dorfman1969economic} and biomedical engineering \cite{cherruault2012mathematical}. The central approach of this theory is the Hamilton-Jacobi-Bellman (HJB) equation, which characterizes the optimal cost-to-go and serves as the foundation for both analytical solutions and numerical approximations in modern control design. Depending on the finiteness of the horizon, the value function satisfies the finite-horizon (FH) or steady-state (infinite-horizon) HJB, characterizing a necessary condition of optimality. The major difference is that the FH-HJB is a time-dependent nonlinear equation involving three unknown terms. There are several extensions of the HJB equation (HJBE) under different settings. Wang et al.~\cite{wang2020reinforcement} introduced the HJBE for an explorative reward function. Wilzer et al. ~\cite {wiltzer2022distributional} proposed the distributional HJBE.

For general nonlinear OCP, the HJBE do not always have a closed-form explicit solution. For the special case of a continuous-time (CT) linear quadratic regulator (LQR), the HJBE can be reduced to an algebraic Riccati equation (ARE), as a matrix-valued equation. For complex scenarios and nonlinear systems, iterative LQR \cite{chen2019autonomous} was designed based on local linearization. For those system which drift term can be decomposed into semi-linear form, state-dependent Riccati equation (SDRE) \cite{ccimen2008state} can provide a reliable suboptimal feedback control maintaining computational efficiency and stabilizing nonlinear dynamics. Several deep neural network (DNN) based nonlinear control algorithms have emerged in recent years. Lin et al. \cite{lin2022policy} exploited FH approximate dynamic programming (FH-ADP)  to solve CT OCP by deriving an explicit time derivative. Nakamura-Zimmerer et al. \cite{nakamura2021adaptive} designed a pure data-driven adaptive framework to solve high-dimensional HJB by utilizing the Pontryagin Minimum Principle (PMP). Albi et al. \cite{albi2021gradient} proposed a gradient-augmented supervised learning-based algorithm to deal with optimal feedback laws using SDRE.  Lee et al. \cite{lee2025hamilton} utilized deep operator learning to accelerate Hamilton-Jacobi-based policy iteration.

The development of scientific machine learning has led to numerous methods that effectively address complex problems. Physics-Informed Neural Networks (PINNs) were introduced by Raissi et al. \cite{raissi2019physics} for solving PDEs, incorporating physical laws as optimization constraints to enable unsupervised training. As a mesh-free technique, PINNs offer flexibility in handling complex boundaries and solving high-dimensional problems at arbitrary collocation points. Subsequent enhancements have focused on improving both accuracy and efficiency for the challenging benchmark problems \cite{wang2022and, anagnostopoulos2024residual, jiao2026gspinn, jiao2024leveraging}. There has been a great interest in the Physics-Informed optimal control strategy. Fotiadis et al. \cite{fotiadis2023physics} regarded the static HJBE as a steady state of an FH-OCP, which can be solved by PINNs as a time-dependent HJBE. Shilova et al. \cite{shilova2024learning} proposed an $\varepsilon$-scheduling iteration to solve the viscosity HJBE by utilizing PINNs. Li et al. proposed an adaptive weighted PINN to solve Euler-Lagrange systems in optimal control \cite{li2026aw}. Banerjee et al. present a survey on physics-informed reinforcement learning (RL) \cite{banerjee2025survey} to clarify how physical priors are integrated in RL. Neural scalar fields also connect PDE-based learning with motion
planning. NTFields~\cite{ni2023ntfields} learns an arrival-time field governed by the Eikonal equation and constructs paths by following its gradients.

In this work, we propose a novel physics-informed method for solving high-dimensional infinite-horizon quadratic cost regulation problems. A neural network parametrizing the value function is trained by integrating the HJBE residual and a control policy label mismatch, using SDRE as a reference. To enhance accuracy and stability, the network architecture enforces zero boundary conditions and positivity. We test the method on high-dimensional nonlinear examples: the 40D Cucker––Smale model and the 51D Allen-–Cahn equation. Numerically, the learned value predictions closely track the
finite-horizon cost-to-go of the SDRE reference controller along the tested trajectories, providing evidence of value--cost consistency. Moreover, the NN controller achieves about 20--30$\times$ faster execution than SDRE. Theoretically, using the viscosity solution of HJB as an intermediate, we identify factors affecting value approximation error, including residual and viscous errors. To further characterize policy error, we derive a Liouville-type inequality for a weighted energy, showing that the weighted policy error is bounded by the weighted residual. The novelty of this paper is summarized as follows.
\begin{itemize}
\item Loss function consists of HJB residual and reference suboptimal control policy, without any requirement for reference value function.
\item This approach can handle high dimensional problems with reliable accuracy, effectively overcoming the heavy computational complexity of traditional SDRE.
\item Considering the mismatch between the approximated value function and cost-to-go, our method achieves lower mismatch than traditional SDRE, indicating improved value--cost consistency. 
\item We provide rigorous error analysis for the approximated value function by NN and control policy.
\end{itemize}

The remainder of the paper is organized as follows. Section \ref{sec:preliminary} presents basic concepts on optimal control, the infinite horizon HJBE. In section \ref{sec:method}, we shall introduce a detailed methodology. Section \ref{sec:analysis} provides a theoretical guarantee about the approximation error of the value function and control policy. Section \ref{sec:experiment} will show numerical examples to demonstrate the effectiveness and accuracy of the proposed algorithm. We conclude the whole paper in section \ref{sec:conclusion}. Appendix can be found to review the detailed proofs and supplementary materials.

\section{Preliminaries}\label{sec:preliminary}

In this section, we introduce the basic formulations for infinite-horizon optimal control.
Our goal is to minimize the following cost function, which depends on the control input $u(t)$ for $t\in[0,\infty)$:
\begin{equation}\label{eq:J}
\begin{aligned}
\min_{u(\cdot)\in\mathbb{U}} J(u(\cdot),x_0) :=& \int_{0}^{\infty} e^{-\lambda s}\,\ell(x(s),u(s))\,{\rm d}s,
\end{aligned}
\end{equation}
where state dynamics $\frac{{\rm d}x}{{\rm d}t} = f(x(t)) + B(x(t))u(t), x(0)=x_0$. Running cost $\ell(x,u):=x^{\top}Qx + u^{\top}Ru$. The state vector $x(t)=[x_1(t),\dots,x_n(t)]^{\top}$ belongs to a bounded domain $\Omega\subset\mathbb{R}^n$. The control signal $u(\cdot)$ takes values in the control space $\mathbb{U}\subset\mathbb{R}^m$, and the admissible control set $\mathbb{U}$ consists of measurable functions mapping $[0,\infty)$ to $\mathbb{U}$.
The exponential discount factor $e^{-\lambda s}$ with $\lambda\ge 0$ ensures convergence of the infinite-horizon integral.
The matrices $Q\in\mathbb{R}^{n\times n}$ and $R\in\mathbb{R}^{m\times m}$ are symmetric positive definite.
The drift dynamics $f:\mathbb{R}^n\to\mathbb{R}^n$ and the control coefficient $B:\mathbb{R}^n\to\mathbb{R}^{n\times m}$ are assumed differentiable.

The value function is defined as $V(x) := \min_{u(\cdot)\in\mathbb{U}} J(u(\cdot),x)$ and the Hamiltonian operator is given by
\begin{equation}\label{Eq:Hamiltonian}
\mathcal{H}(x,u,\nabla V) := \nabla V(x)^{\top}\bigl(f(x)+B(x)u\bigr) + \ell(x,u).
\end{equation}

When the control space is unconstrained, i.e., $\mathbb{U}=\mathbb{R}^m$, the optimal control policy admits a closed-form expression:
\begin{equation}\label{OC_u}
u(x,\nabla V) = -\frac{1}{2}R^{-1}B(x)^{\top}\nabla V(x).
\end{equation}
The discounted HJBE reads
\[
\lambda V(x) - \min_{u\in\mathbb{U}} \mathcal{H}(x,u,\nabla V) = 0,
\]

and substituting the optimal policy~\eqref{OC_u} yields the HJBE satisfied by the value function:
\begin{equation}\label{HJB_V}
\left\{
\begin{aligned}
&\lambda V(x) - x^{\top}Qx - \nabla V(x)^{\top} f(x) + \frac{1}{4}\nabla V(x)^{\top} B(x)R^{-1}B(x)^{\top}\nabla V(x) = 0, && x\in\Omega,\\
&V(0) = 0.
\end{aligned}
\right.
\end{equation}
When the discount factor vanishes, $\lambda=0$, equation~\eqref{HJB_V} reduces to the static HJBE involving only the gradient $\nabla V$.
Once $V(x)$ is determined, the optimal feedback law is given by~\eqref{OC_u}.

An alternative description of the system dynamics is the semi-linear form
\begin{equation}\label{eq:dyn}
\frac{{\rm d}x}{{\rm d} t} = A(x)x + B(x)u(t),
\end{equation}
where $A(x)$ is a state-dependent matrix.
This structure enables the state-dependent Riccati equation (SDRE) approach, which provides a suboptimal feedback law by solving the algebraic Riccati equation at each state $x$:
\begin{equation}\label{state_ARE}
A(x)^{\top}P(x) + P(x)A(x) - P(x)B(x)R^{-1}B(x)^{\top}P(x) + Q - \lambda P(x) = 0.
\end{equation}
Analogously to the linear quadratic regulator, the SDRE-based approximation of the value function is constructed as
\begin{equation}\label{sdre_V}
V_S(x) = x^{\top}P(x)x,
\end{equation}
where $P(x)$ satisfies~\eqref{state_ARE} and is positive semi-definite.
The corresponding exact feedback control is
\begin{equation}\label{SDRE-control}
u_{V_S}(x)
:=-\frac12 R^{-1}B(x)^\top\nabla V_S(x)
=-R^{-1}B(x)^\top
\left[P(x)x+\frac12\varphi(x)\right],
\end{equation}

with the components $\varphi_i(x) = x^{\top}\frac{\partial P(x)}{\partial x_i}x, 1\le i\le n$. The standard SDRE feedback used to generate the control labels is
\begin{equation}\label{SDRE-control-approx}
u_S(x):=-R^{-1}B(x)^\top P(x)x.
\end{equation}
This feedback omits the derivatives of $P(x)$.

\section{Methodology}\label{sec:method}
We adopt a PINN approach for high-dimensional nonlinear optimal control problems.
The value function $V(x)$ is approximated by a neural network.
This approach is motivated by several key advantages:
First, learning a scalar-valued $V(x)$ aligns with dynamic programming principles and the feedback policy is derived from the value gradient;
Second, approximating a scalar field mitigates the curse of dimensionality more effectively than learning a vector field $\nabla V(x)$, improving sample efficiency;
Third, the gradient $\nabla V(x)$ required for the control law is computed accurately and efficiently via automatic differentiation, enabling fast online evaluation;
Finally, the framework incorporates system physics by directly embedding the HJB residual into the loss function.

The subsequent subsections detail the key components of this method.
\subsection{Data Collection}

We construct the dataset $\mathcal{D}$ by simulating short trajectories under the SDRE controller. First, $N_{\mathrm{traj}}$ initial states $\{x_0^{(j)}\}_{j=1}^{N_{\mathrm{traj}}}$ are drawn independently from a chosen distribution over the region of interest.
For each initial state, we iterate the following procedure for $N_T$ steps:
at step $k$ ($0 \le k < N_T$), the state is $x_k^{(j)}$; we compute the SDRE-based control $u_k^{(j)} = u_S(x_k^{(j)})$ via~\eqref{SDRE-control-approx}; then the system is advanced by a time increment $\Delta t$ to obtain the next state $x_{k+1}^{(j)}$, using the dynamics~\eqref{eq:dyn} under control $u_k^{(j)}$.
The integration is performed with a standard numerical solver (e.g., Runge–Kutta RK45).
In the end, we collect all state-control pairs:
\[
\mathcal{D} = \bigl\{(x_i, u_i)\bigr\}_{i=1}^{N},\qquad
N = N_{\mathrm{traj}} \times N_T,
\]
where the index $i$ runs over all trajectories and time steps. Under such a sampling strategy, the training process will be implemented more on the state region visited by reference trajectories.

\subsection{Physics-informed optimal feedback algorithm}

Let $V_\theta(x)$ denote the neural network approximation of the value function with parameters $\theta = \{W_i, b_i\}$ including weights and bias.
The optimal control policy is approximated from the gradient of $V_\theta$ as $u(x,\nabla V_{\theta})$, where $u$ is defined in \eqref{OC_u}.
A key limitation of the classical SDRE approach is the heavy computational burden of solving a Riccati equation for $P(x)$ at every state encountered along a trajectory.
Our framework alleviates this burden by using the SDRE control only as a soft label during training, while simultaneously enforcing the HJBE as a physical constraint.
The loss function comprises two terms: the HJB term and the control policy term.

\paragraph{HJB term.}
This term enforces the HJB optimality condition by penalizing its residual:
\begin{equation}\label{eq:Lr}
\mathcal{L}_{\rm HJB}(\theta) = \frac{1}{N} \sum_{i=1}^{N} \bigl| \lambda V_\theta(x_i) - \mathcal{H}\bigl(x_i, u_{\theta}(x_i), \nabla V_\theta(x_i)\bigr) \bigr|^2,
\end{equation}
where $\mathcal{H}$ is the Hamiltonian defined in \eqref{Eq:Hamiltonian}, $u_{\theta}(x_i)= u(x_i,\nabla V_{\theta}(x_i))$, and $u$ is defined in \eqref{OC_u}. 

\paragraph{Control policy term.}
This term encourages the network output to match the near-optimal SDRE control signals from the dataset:
\begin{equation}\label{eq:Lc}
\mathcal{L}_{c}(\theta) = \frac{1}{N} \sum_{i=1}^{N} \bigl| u_S(x_i) - u_\theta(x_i) \bigr|^2,
\end{equation}
where $u_S$ is given by~\eqref{SDRE-control-approx}.

The total training loss is a weighted sum of the two terms:
\begin{equation}\label{Loss_eq}
\mathcal{L}(\theta) = \mathcal{L}_{\rm HJB}(\theta) + \gamma \mathcal{L}_{c}(\theta),
\end{equation}
where $\gamma > 0$ is a tuning parameter that balances the influence of the HJB optimality condition and the data fidelity provided by the SDRE labels.

\paragraph{Network parameterization.}
To improve numerical conditioning, we normalize the inputs using
full-pool state statistics.
Let $\mu\in\mathbb{R}^d$ and $\sigma\in\mathbb{R}^d$ denote the
coordinate-wise empirical mean and standard deviation computed from
the collected state dataset before the pointwise
training/held-out split. The held-out states contribute to these preprocessing statistics;
closed-loop performance is evaluated using separate test initial
conditions. We define the normalized state as $\tilde{x} = \frac{x - \mu}{\sigma},
\tilde{x}_0 = \frac{-\mu}{\sigma}.$
We then construct $z_\theta(x) = \mathcal{NN}_\theta(\tilde{x}) - \mathcal{NN}_\theta(\tilde{x}_0)$, where $\mathcal{NN}_\theta(\cdot)$ is a standard feedforward neural network parameterized by $\theta$.
The value function is given by
\[
V_\theta(x) = z_\theta(x) + \bigl( \mathrm{ReLU}(z_\theta(x)) - z_\theta(x) \bigr)_{\mathrm{sg}},
\]
where $\mathrm{sg}(\cdot)$ denotes the stop-gradient operator.
By construction, $V_\theta(0)=0$ and $V_\theta(x) \ge 0$ for all $x$ since $V_\theta(x) = \mathrm{ReLU}(z_\theta(x))\ge 0$ and $\nabla_\theta V_\theta (x) = \nabla_\theta z_\theta (x)$. Two operations, i.e.,  ReLU and the stop-gradient operation, ensure positivity and invariant gradient calculation, respectively. In practice, the SDRE-guided training reliably drives $z_\theta$ into the non-negative region, making the stop-gradient a harmless safeguard.



\paragraph{Training procedure.}
The network is trained by minimizing the total loss~\eqref{Loss_eq} over the dataset $\mathcal{D}$ using a stochastic gradient-based optimizer.
The training procedure is summarized in Algorithm~\ref{alg:pinn_sdre}.

\begin{algorithm}[H]\small
\caption{SDRE-guided PINN for optimal control}
\label{alg:pinn_sdre}
\KwIn{domain $\Omega\subset\mathbb{R}^n$, system parameters $f, A, B, Q, R$, discount factor $\lambda$, weight $\gamma$, dataset $\mathcal{D} = \{(x_i, u_i)\}_{i=1}^N$.
}

Initialize neural network parameters $\theta$.

\For{epoch $k = 1$ to $\text{MAX\_ITER}$}{
    \For{mini-batch $S \subset \mathcal{D}$}{
        Let $S_x = \{x_i : (x_i, u_i) \in S\}$.

        Compute HJB term $\mathcal{L}_{\rm HJB}$ on $S_x$ via~\eqref{eq:Lr}.
        
        Compute control policy term $\mathcal{L}_c$ on $S_x$ via~\eqref{eq:Lc}.

        $\mathcal{L} \leftarrow \mathcal{L}_{\rm HJB} + \gamma \mathcal{L}_c$.
        
        Update $\theta \leftarrow \theta - \eta \nabla_\theta \mathcal{L}$.
    }
}

{\bf Return} Approximate value function $V_{\theta}$.
\end{algorithm}



\section{Error analysis}\label{sec:analysis}
In this section, we shall focus on the estimate of error between the NN solution $V_\theta$ and the exact value function $V$. From a PDE perspective, it is relatively hard to directly consider $V$, which might be non-smooth and have multiple solutions. A good choice here we considered is to introduce a viscosity solution $V^\vare$ as an intermediate. In the first part, we show that under certain bounded residual conditions, the error itself $|V-V_\theta|$ can be arbitrarily close to zero. In the second part, the gradient error $|\nabla V - \nabla V_\theta|$ will be considered by applying PDE-related estimation. This term quantifies a certain metric to compare the NN approximated policy and the optimal policy.

\subsection{Generalization error of value function}
For a smooth function $U$, define the viscous HJB residual by
\begin{equation}\label{eq:viscous-residual}
r_\varepsilon[U](x)
:=\lambda U(x)
-\min_{u\in\mathbb{U}}\mathcal{H}(x,u,\nabla U)
-\varepsilon\Delta U(x).
\end{equation}

For the SDRE value function approximation \eqref{sdre_V}, there is an explicit expression of residual.
\begin{equation}\label{residual_VS}
\begin{aligned}
r_\varepsilon[V_S](x)
={}&-\varphi(x)^\top
[A(x)-B(x)R^{-1}B(x)^\top P(x)]x
+\frac14\varphi(x)^\top
B(x)R^{-1}B(x)^\top\varphi(x)\\
&-\varepsilon\left[
4(\nabla\cdot P(x))^\top x
+2\operatorname{tr}(P(x))
+x^\top\Delta P(x)x
\right].
\end{aligned}
\end{equation}
This expression follows by differentiating~\eqref{sdre_V}
and using~\eqref{state_ARE}. 
The notations $[\varphi(x)]_i = x^\top \partial_{x_i} P(x)x$, $(\nabla\cdot P(x))_j := \sum_{k=1}^n \frac{\partial P_{kj}(x)}{\partial x_k}$ and $(\Delta P(x))_{ij} := \sum_{k=1}^n \frac{\partial^2 P(x)}{\partial x_k^2}$.

For neural network approximation $V_\theta$, similarly we have corresponding residual $r_\vare[V_\theta](x)$. Rewrite the equation, it can be seen that $V_\theta$ satisfies the viscousity solution of HJB with an extra residual,
\begin{equation}\label{Eq:Vtheta}
\lambda V_\theta(x)
-\min_{u\in\mathbb{U}}
\left\{
\mathcal{H}(x,u,\nabla V_\theta)
+r_\varepsilon[V_\theta](x)
\right\}
-\varepsilon\Delta V_\theta(x)=0.
\end{equation}
Notice that this is a second-order equation with the term $-\vare \Delta V_\theta$, which can be interpreted in a stochastic optimal control framework hereafter. Let us define a modified running cost $\tilde{\ell}_\varepsilon(x,u)
=x^\top Qx+u^\top Ru+r_\varepsilon[V_\theta](x)$. With this choice, we get the second-order modified viscousity form of the HJBE accordingly.
\begin{equation}\label{Eq:modified-viscous-HJB}
\lambda V_\theta(x)
-\min_{u\in\mathbb{U}}
\left\{
\nabla V_\theta(x)^\top(f(x)+B(x)u)
+\tilde{\ell}_\varepsilon(x,u)
\right\}
-\varepsilon\Delta V_\theta(x)=0.
\end{equation}
The solution of the above equation corresponds to the following stochastic optimal control problem with discount factor:
\[
\left\{
\begin{aligned}
\min_{u}\tilde J_\vare(x; u) =& \mathbb{E}\bigg[\int_0^\infty e^{-\lambda t}\tilde\ell_\vare(X_t, u_t)dt\bigg],\\
\text{s.t.} \quad dX_t =& [f(X_t)+B(X_t)u_t]dt + \sqrt{2\vare}dW_t.
\end{aligned}
\right.
\]
In the following, we show that $V_\theta(x)$ can be equivalently described by the stochastic dynamical programming principle (SDPP). The detailed proof can be found in section \ref{pf:Lem:SDPP}.

\begin{lemma}\label{Lem:SDPP}
Let $V_\theta\in C^2(\mathbb{R}^n)$ denote the neural-network
value approximation that satisfies the modified viscous HJBE \eqref{Eq:modified-viscous-HJB}. Then we have the Stochastic Dynamical Programming Principle (SDPP).
\[
    V_\theta(x) = \inf_{u\in\mathbb{U}}\mathbb{E}\left[\int_0^T e^{-\lambda s}\tilde{\ell}_\vare(X_s, u_s)ds + e^{-\lambda T}V_\theta(X_T)\right].
\]
\end{lemma}

In the following result, we shall provide an upper bound for the error $|V_\theta - V^\vare|$. Detailed proof can be found in the section \ref{pf:prop:VS-Ve}.
\begin{theorem}\label{prop:VS-Ve}
Let $u^*$ be an optimal control for $V^\varepsilon$, and let $\widetilde u$ be an optimal control for the modified problem associated with $V_\theta$. Denote by $X_t^*$ and $\widetilde X_t$ the corresponding state processes. Assume furthermore that
$\lim_{T\to\infty}\mathbb E_x\!\left[e^{-\lambda T}\bigl|V_\theta(X_T^*)-V^\varepsilon(X_T^*)\bigr|\right]=0$ and $\lim_{T\to\infty}\mathbb E_x\!\left[e^{-\lambda T}\bigl|V_\theta(\widetilde X_T)-V^\varepsilon(\widetilde X_T)\bigr|\right]=0$.  Then
\[
    |V_\theta(x)-V^\varepsilon(x)|
\le
\max\left\{
\mathbb E_x\!\int_0^\infty e^{-\lambda t}|r_\varepsilon(X_t^*)|\,dt,\,
\mathbb E_x\!\int_0^\infty e^{-\lambda t}|r_\varepsilon(\widetilde X_t)|\,dt
\right\}.
\]
In particular, if $\lambda>0$ and $r_\varepsilon\in L^\infty$, then $|V_\theta(x)-V^\varepsilon(x)|
\le
\frac{1}{\lambda}\|r_\varepsilon\|_{L^\infty}$.
\end{theorem}

\begin{corollary}
Let $V$ be the optimal value function defined by the standard HJBE. We have the following

\noindent{(i)} If $\lambda>0$ and $r_\vare\in L^\infty$, the value error $|V_\theta(x)-V(x)|\le
\frac{1}{\lambda}\|r_\varepsilon\|_{L^\infty} + C\sqrt{\vare}$;

\noindent{(ii)} If $\lambda = 0$, then $|V_\theta(x)-V(x)|
\le
\max\left\{
\mathbb E_x\!\int_0^\infty |r_\varepsilon(X_t^*)|\,dt,\,
\mathbb E_x\!\int_0^\infty |r_\varepsilon(\widetilde X_t)|\,dt
\right\} + C\sqrt{\vare}$.

\end{corollary}

The proof can be easily obtained by combining the result of Theorem \ref{prop:VS-Ve} and the well-known estimate \eqref{Eq:V-Ve} together. On the right-hand side, the first term reflects the mismatch between NN $V_\theta$ for viscosity solution $V^\vare$. The second term gives the convergence rate when considering a vanishing viscous sequence $V^\vare \to V$ as $\vare \to 0$.

\subsection{Gradient estimate and policy error}
In deep learning-based optimal control, quantifying how well a neural network approximates the control policy is critical. Since the original HJBE is difficult to handle directly, we instead analyze its viscosity form, which facilitates energy estimates.

The key approximation error is the gradient $|\nabla V^\varepsilon - \nabla V_\theta|$. Let $w(x) = V^\varepsilon(x) - V_\theta(x)$ denote the error between the neural network and the exact viscosity solution. We define $E(r)$ as a weighted $L^2$-norm of the gradient error over a domain of radius $r$, with weight determined by $f(x)$ and $V^\varepsilon$, and derive a Liouville-type differential inequality for $E(r)$.

Hereafter, assuming that the weighted residual is sufficiently small and that the weighted boundary measure is uniformly bounded, a contradiction argument shows that $E(r)$ cannot exceed the prescribed tolerance $\delta$. Detailed proof can be found in the section \ref{pf:Thm-error-gradient-estimate}.

\begin{theorem}\label{Thm-error-gradient-estimate}
Assume that the control matrix $B(x)$ has full rank and  $B(x)R^{-1}B(x)^\top=I$ is an identity matrix. Assume the drift function has gradient form $f=\nabla \phi$. Denote $\Phi(x):=\frac12 V^\varepsilon(x)-\phi(x)$ and $\rho(x):=\exp\left(-\frac{\Phi(x)}{\varepsilon}\right)$. Let total residual error $R_\varepsilon(x):=r_\varepsilon[V_\theta](x)+\lambda w(x)$. Denote $M(r) := \int_{B_r} \rho$ and $H(r) := \int_{B_r} \rho |R_\vare|$. Assume 

(i) $M'(r)=\int_{\partial B_r}\rho \le C$;\quad 
(ii) $\left\|R^{-1}B(x)^{\top} \right\| \le K_r, x\in B_r$; \quad
(iii) $\int_{B_r}\rho |\nabla (V^\vare - V)|^2\le \eta_\vare(r)$. 

If the global weighted residual error satisfies $H_\infty:=\int_{\mathbb{R}^n}\rho |R_\varepsilon|
\le \frac{\delta}{8}$, then for every $r>0$,
We shall have the policy estimate:
\[
\int_{B_r} \rho |u_\theta - u|^2 \le \frac{K_r^2}{2}(\delta + \eta_\vare(r)).
\]
\end{theorem}

Let us explain the connection between the training process of NN and the above upper bound. 

Let $S_N=\{x_i\}_{i=1}^N$ be the training points with quadrature weights $\{\omega_i\}_{i=1}^N$.  Define the empirical weighted source loss by $\mathcal E_T^R
:=
\sum_{i=1}^N
\omega_i
\rho(x_i)
\left|
r[V_\theta](x_i)
\right|$, where $r_\vare[V_\theta] = r[V_\theta] -\varepsilon\Delta V_\theta$. Assume that the quadrature rule satisfies $\left|
\int_{\mathbb R^n}
\rho |r[V_\theta](x)|
-
\sum_{i=1}^N
\omega_i
\rho(x_i)
|r[V_\theta](x_i)|
\right|
\le
C_{\mathrm{quad}}N^{-\alpha}.$
Then we conclude the upper bound for the weighted residual:
\[
H_\infty
\le
\underbrace{\mathcal E_T^R}_{\text{training error}}
+
\underbrace{C_{\mathrm{quad}}N^{-\alpha}}_{\text{quadrature error}} + \underbrace{\mathcal{O}(\vare)}_{\text{viscosity error}} + \underbrace{|\lambda|\cdot\|w\|_{L^1_\rho(\R^n)}}_{\text{value error}}.
\]
This is the direct analogue of the PINNs generalization estimate: the global weighted source error $H_\infty$ is upper bounded by the empirical training source loss, the quadrature error, viscosity error, and value error.

The full gradient estimate $\nabla V$ is relatively difficult to obtain due to the quadratic term in HJBE. Instead, if estimating the directional gradient $B^\top\nabla V$, the control policy could be directly upper bounded without several conditional assumptions in Theorem \ref{Thm-error-gradient-estimate}. A detailed proof appears in Section \ref{pf:Thm-policy-estimate}.
\begin{theorem}\label{Thm-policy-estimate}
Let $V^\varepsilon$ be a classical solution of the viscous HJB equation, and let $V_\theta$ be the neural network approximation, with
$V^\varepsilon,V_\theta\in C^2(\Omega)$.
Assume $\varepsilon>0$, $\lambda\geq 0$, and
$R\succ 0$. Define $G(x)=B(x)R^{-1}B(x)^\top,
w(x)=V^\varepsilon(x)-V_\theta(x),$
and $b_0(x)=-f(x)+\frac12G(x)\nabla V^\varepsilon(x)\in W^{1,\infty}_{\mathrm{loc}}(\Omega)$. Let $u_\theta = u(x, \nabla V(x))$ and $u^\varepsilon:=u(x, \nabla V^\varepsilon(x))$ defined in \eqref{OC_u} For any $B_r\Subset B_{r_1}\Subset\Omega$, there exists
a constant $C_{\varepsilon,r,r_1}>0$, independent of
$V_\theta$ such that if the vanishing-viscosity policy error
satisfies $\|u^\varepsilon-u\|_{L^2(B_r)}^2
\leq \eta_\varepsilon(r),$
then
\begin{equation}
\label{eq:total_policy_error}
\|u_\theta-u\|_{L^2(B_r)}^2
\leq
\frac{2}{\lambda_{\min}(R)}
\Big(
\|r_\varepsilon[V_\theta]\|_{L^1(B_{r_1})}+C_{\varepsilon,r,r_1}
\|V^\varepsilon-V_\theta\|_{L^1(B_{r_1})}
\Big)
+2\eta_\varepsilon(r).
\end{equation}
\end{theorem}

\section{Numerical Experiments}\label{sec:experiment}
For the main-text experiments, we set $\gamma=1$ in the loss~\eqref{Loss_eq}. 
For each benchmark, the value network has three hidden
layers of 256 $\tanh$ units and is trained with Adam. Initialization and split seeds are fixed at 42. The collected state--control samples are split pointwise
into training and held-out sets at ratios of 80/20 for
Cucker--Smale and 90/10 for Allen--Cahn.
Cucker--Smale uses 4000 epochs, batch size 128, and an initial
learning rate of $2\times10^{-4}$; Allen--Cahn uses 3000 epochs,
batch size 4096, and an initial learning rate of
$2\times10^{-3}$. All simulations are implemented in \textbf{Pytorch} and run on a laptop with an AMD Ryzen 9 8945HX (16 cores, 2.50 GHz base, 3.50 GHz boost) and an NVIDIA GeForce RTX 5060 Laptop GPU (8 GB dedicated memory). \footnote{Code and the datasets and trained models used in
the main-text experiments are available at
\url{https://github.com/myspark2026/PI_optimal-control-by-control-only-supervision}.}

For a given trajectory, we denote by \(J_T^*(x(t))\) the cost accumulated by the SDRE controller from state $x(t)$ over the finite horizon $[t, T]$. This quantity serves as a reference in comparisons.


\noindent{\bf Cucker-Smale model.}\; We consider a multi-agent system with $N_a = 20$ agents governed by:
\begin{equation}
    \begin{bmatrix} \dot{y} \\ \dot{v} \end{bmatrix} = 
    \begin{bmatrix}
        0_{N_a} & I_{N_a} \\
        0_{N_a} & A_{N_a}(y)
    \end{bmatrix}
    \begin{bmatrix} y \\ v \end{bmatrix} + 
    \begin{bmatrix} 0_{N_a} \\ I_{N_a} \end{bmatrix} u,
\label{eq:discrete_system_test1}
\end{equation}
where $0_n$ denotes the $n\times n$ zero matrix, and $ [A(y)]_{i,j} = -\dfrac{1}{N_a} \displaystyle\sum_{k \neq i} P(y_i, y_k)$ if $i=j$ otherwise $\dfrac{1}{N_a} P(y_i, y_j)$.
with $P(y_i, y_j) = 1/({1 + \|y_i - y_j\|^2})$.
The cost functional is $ J(y,v,u) = \frac{1}{N_a} \int_0^T \sum_{i=1}^{N_a} \bigl( \|y_i\|^2 + \|v_i\|^2 + \|u_i\|^2 \bigr) \, {\rm d}t$. To generate the training set $\mathcal{D} = \{(x_i, u_i)\}_{i=1}^N$, we sample initial states $x_0 = [y_0^\top,v_0^\top]^\top$ uniformly from $(-1,1)^{40}$. For each of $N_{\text{traj}} = 400$ independent initial conditions, we simulate the closed-loop system under the SDRE controller for $N_T = 50$ steps with a fixed $\Delta t=0.2$; state-control pairs are recorded at every step, yielding $N = 20\,000$ samples in total.

\begin{figure}[t]
\centering
\begin{minipage}[b]{0.04\linewidth}
    \centering
    \rotatebox{90}{\hskip0.5in Position}
\end{minipage}
\;
\begin{minipage}[b]{0.45\linewidth}
    \centering
    NN controller\\
    \includegraphics[width=\linewidth]{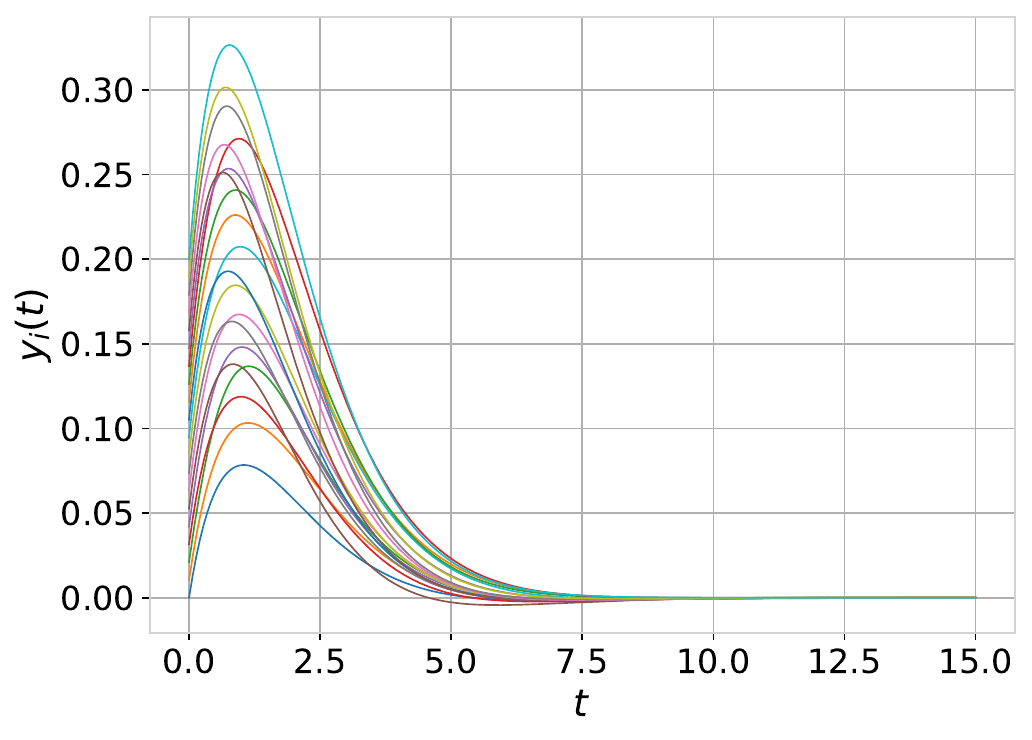}
\end{minipage}
\;
\begin{minipage}[b]{0.45\linewidth}
    \centering
    SDRE controller\\
        \includegraphics[width=\linewidth]{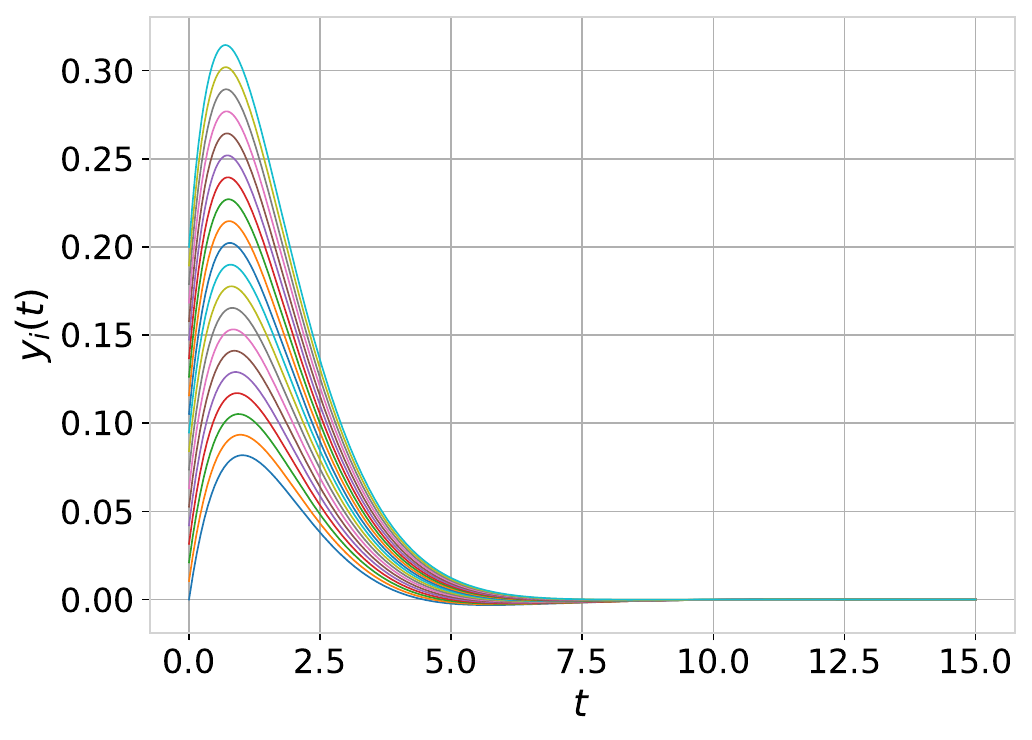}
\end{minipage}

\begin{minipage}[b]{0.04\linewidth}
    \centering
    \rotatebox{90}{\hskip0.5in Velocity}
\end{minipage}
\;
\begin{minipage}[b]{0.45\linewidth}
    \centering
    \includegraphics[width=\linewidth]{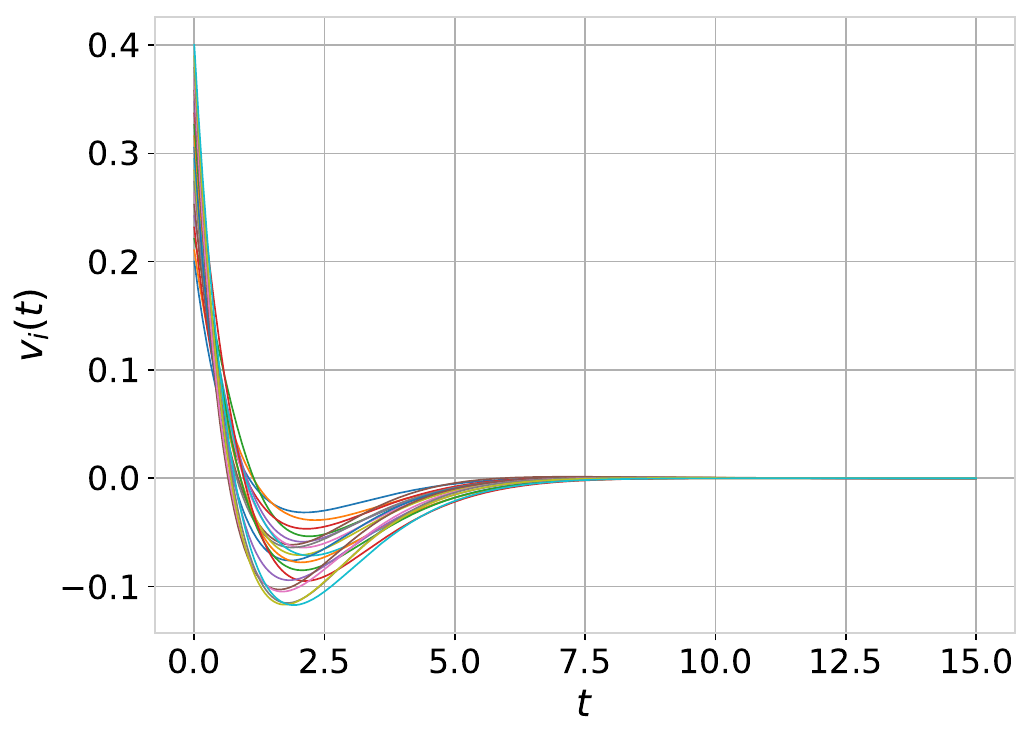}
\end{minipage}
\;
\begin{minipage}[b]{0.45\linewidth}
    \centering
    \includegraphics[width=\linewidth]{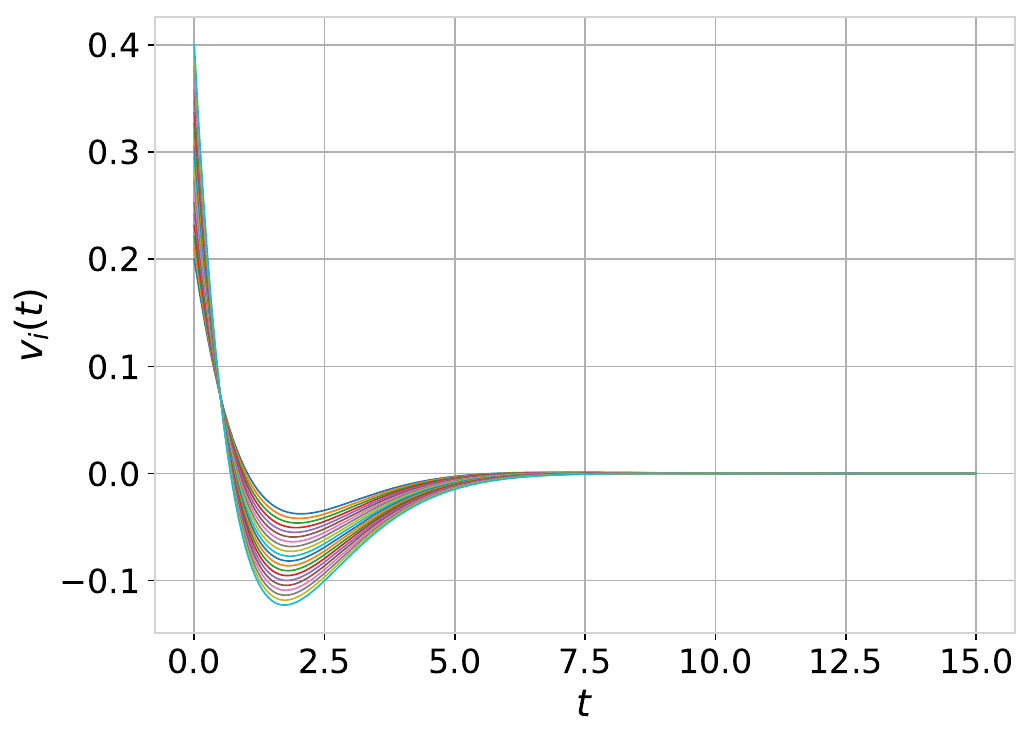}
\end{minipage}

\caption{
Position and velocity trajectories of the Cucker-Smale model under the NN controller (left column) and the SDRE controller (right column).
}
\label{fig:CuckerS_traj}
\end{figure}

\begin{figure}[t]
\centering
\begin{minipage}[b]{0.45\linewidth}
    \centering
    \includegraphics[width=\linewidth]{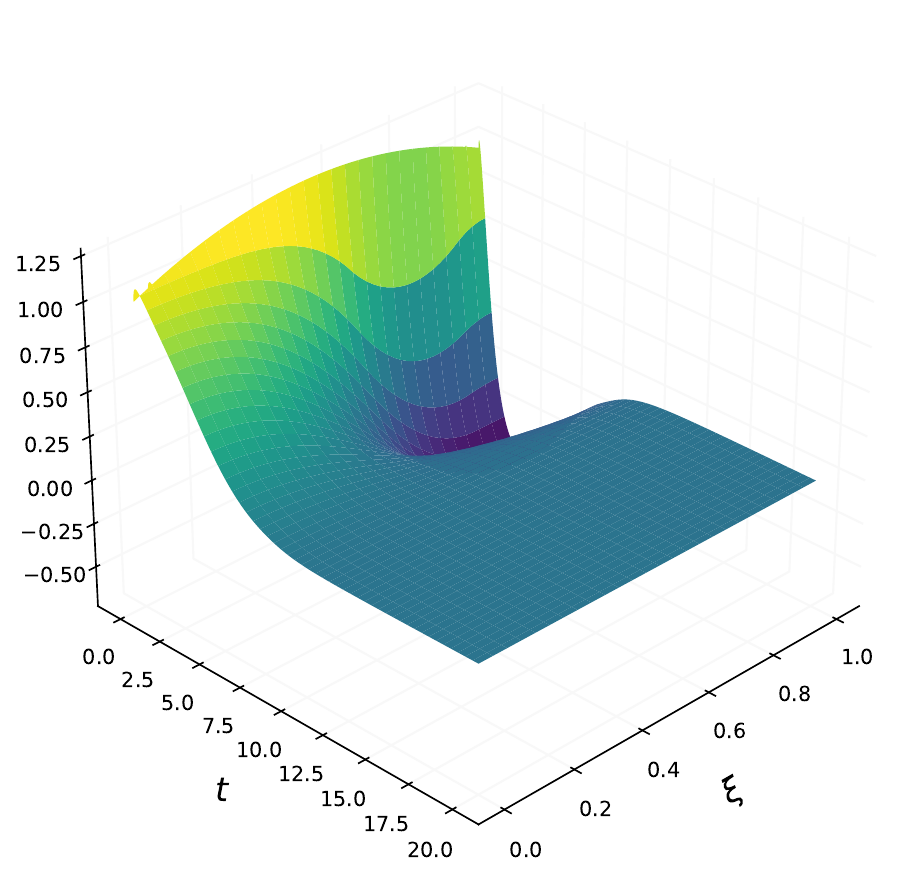}
\end{minipage}
\quad
\begin{minipage}[b]{0.45\linewidth}
    \centering
    \includegraphics[width=\linewidth]{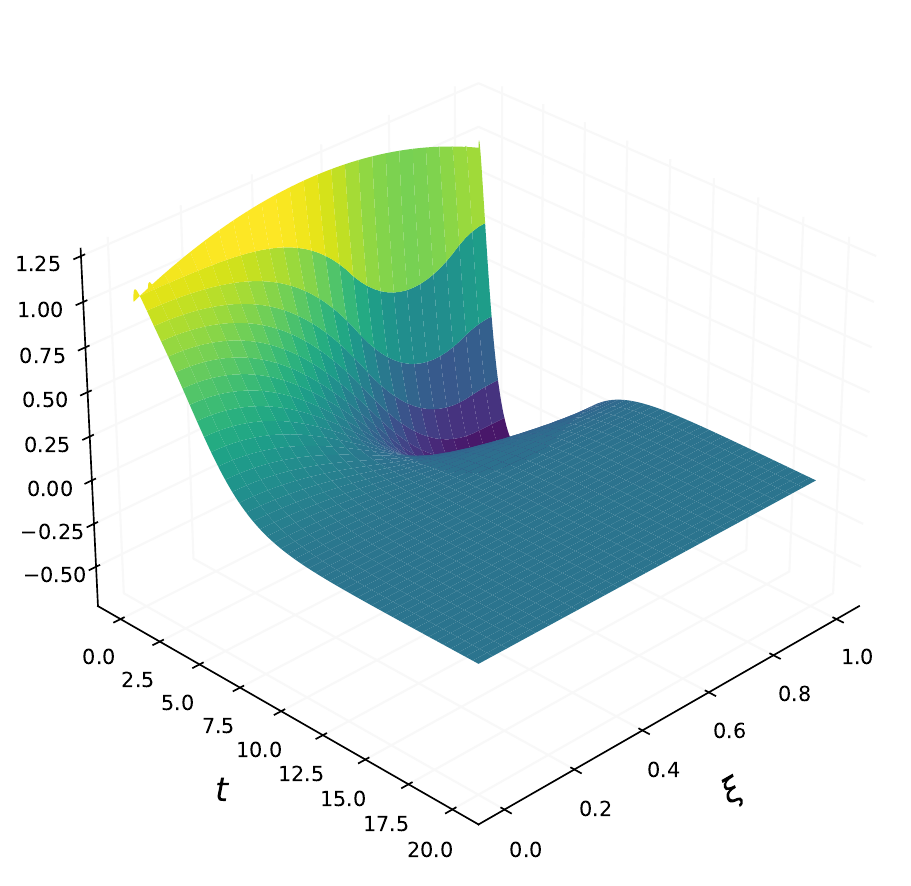}
\end{minipage}
\caption{
Solution of the Allen-Cahn equation under the NN controller (left) and the SDRE controller (right).
}
\label{fig:Allen-Cahn state}
\end{figure}

\begin{figure}[t]
\centering
\begin{minipage}[b]{0.32\linewidth}
    \centering
    \includegraphics[width=\linewidth]{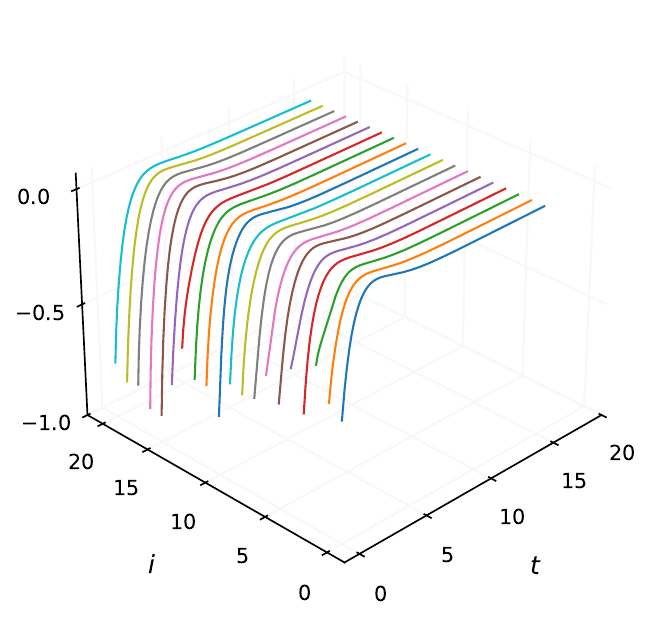}
\end{minipage}
\hfill
\begin{minipage}[b]{0.32\linewidth}
    \centering
    \includegraphics[width=\linewidth]{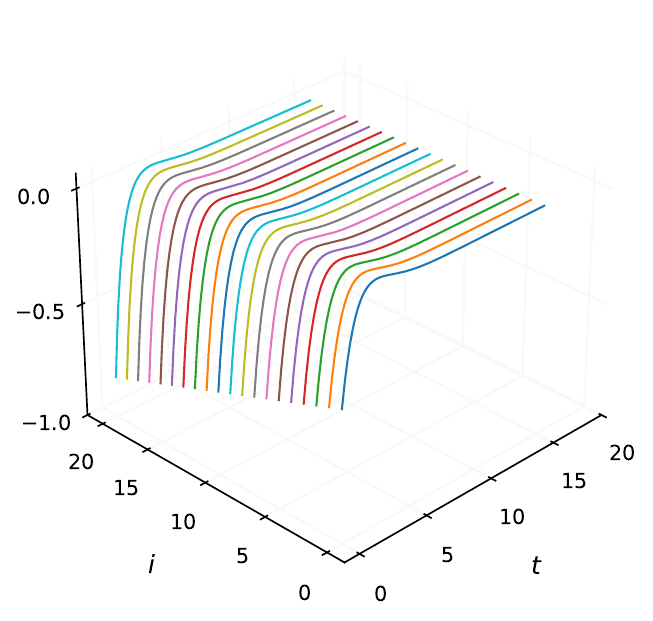}
\end{minipage}
\hfill
\begin{minipage}[b]{0.32\linewidth}
    \centering
    \includegraphics[width=\linewidth]{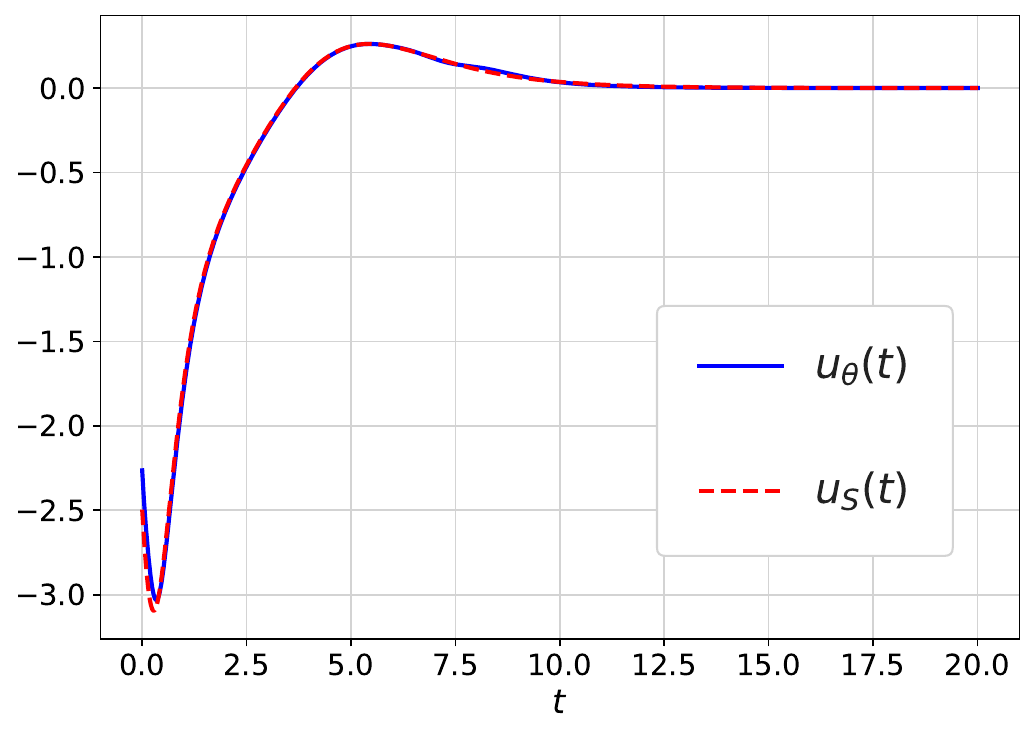}
\end{minipage}
\caption{Control signals. Left: Cucker--Smale under NN controller. Middle: Cucker--Smale under SDRE controller. Right: Allen--Cahn under NN and SDRE controllers.}
\label{fig:control_signal}
\end{figure}

\noindent{\bf Allen--Cahn equation.}\;
We consider the nonlinear Allen--Cahn PDE with Neumann boundary conditions~\cite{albi2021gradient}:
\begin{equation*}
\partial_t x(\xi, t) = 0.1 \partial_{\xi\xi} x + x(1 - x^2) + \chi_\omega(\xi) u(t),
\end{equation*}
where $\xi \in [0,1]$ and $t \in (0,\infty)$.
The control $u : [0,\infty) \to \mathbb{R}$ is scalar, and $\chi_{\omega}(\xi)$ is the indicator function of the interval $\omega = [0.6,0.9]$.
We aim to minimize $\mathcal{J}(u, x) = \int_0^{+\infty} \left( \int_0^1 |x(\xi,t)|^2 d\xi + 0.1 u^2(t) \right) {\rm d}t$. Using finite differences on a uniform grid of $d = 51$ points, the PDE is semi-discretized as
\begin{equation}
\dot{{x}}(t) = {A}({x}) {x}(t)  + {B} u(t),
\quad\mbox{with}\quad
  {A} ({x}) = 0.1 A_0 + I_d -diag({x} \odot {x}),
  \label{eq:discrete_system_test2}
\end{equation}
where $\odot$ is the Hadamard product, $A_0$ is the tridiagonal discrete Laplacian matrix and
\[
{B}_i = \chi_{[0.6,0.9]}(\xi_i) =
\begin{cases}
1, & \xi_i \in [0.6, 0.9], \\
0, & \text{otherwise},
\end{cases}
\qquad i=1,\dots,d.
\]

We generate $N_{\text{traj}} = 50$ training trajectories. The initial profiles are constructed as $x_0(\xi) = s \sum_{k=1}^{5} a_k \cos(2\pi k \xi + \phi_k) + b$, where $a_k \sim \mathcal{U}(-\frac{1}{k}, \frac{1}{k})$, $\phi_k \sim \mathcal{U}(0,2\pi)$, $s \sim \mathcal{U}(0.5,2)$, and $b \sim \mathcal{U}(-1,1)$.
The random scaling factor $s$ and baseline shift $b$ provide diverse initial conditions. For each t initial condition, we simulate $N_T = 200$ steps with $\Delta t=0.01$ under the SDRE controller and collect state–control pairs. This gives $N = 10\,000$ samples.

The Cucker--Smale and Allen--Cahn benchmarks further evaluate the proposed method beyond the structural assumptions of Theorem~\ref{Thm-error-gradient-estimate}. Specifically, $B\in\mathbb{R}^{40\times20}$ and $B\in\mathbb{R}^{51\times1}$, respectively, so neither system satisfies $BR^{-1}B^\top=I$. Importantly, this condition is required only for the weighted policy-error analysis, not for the proposed training algorithm. Moreover, Theorem~\ref{prop:VS-Ve} and its corollary establish value-error bounds under a separate set of assumptions without requiring this structural condition. While the residual conditions underlying these bounds are not directly verified by finite trajectory samples, the numerical results provide complementary empirical evidence of the method's effectiveness in the tested operating regions.

\begin{figure}[!hbt]
\centering
\includegraphics[width=.38\linewidth]{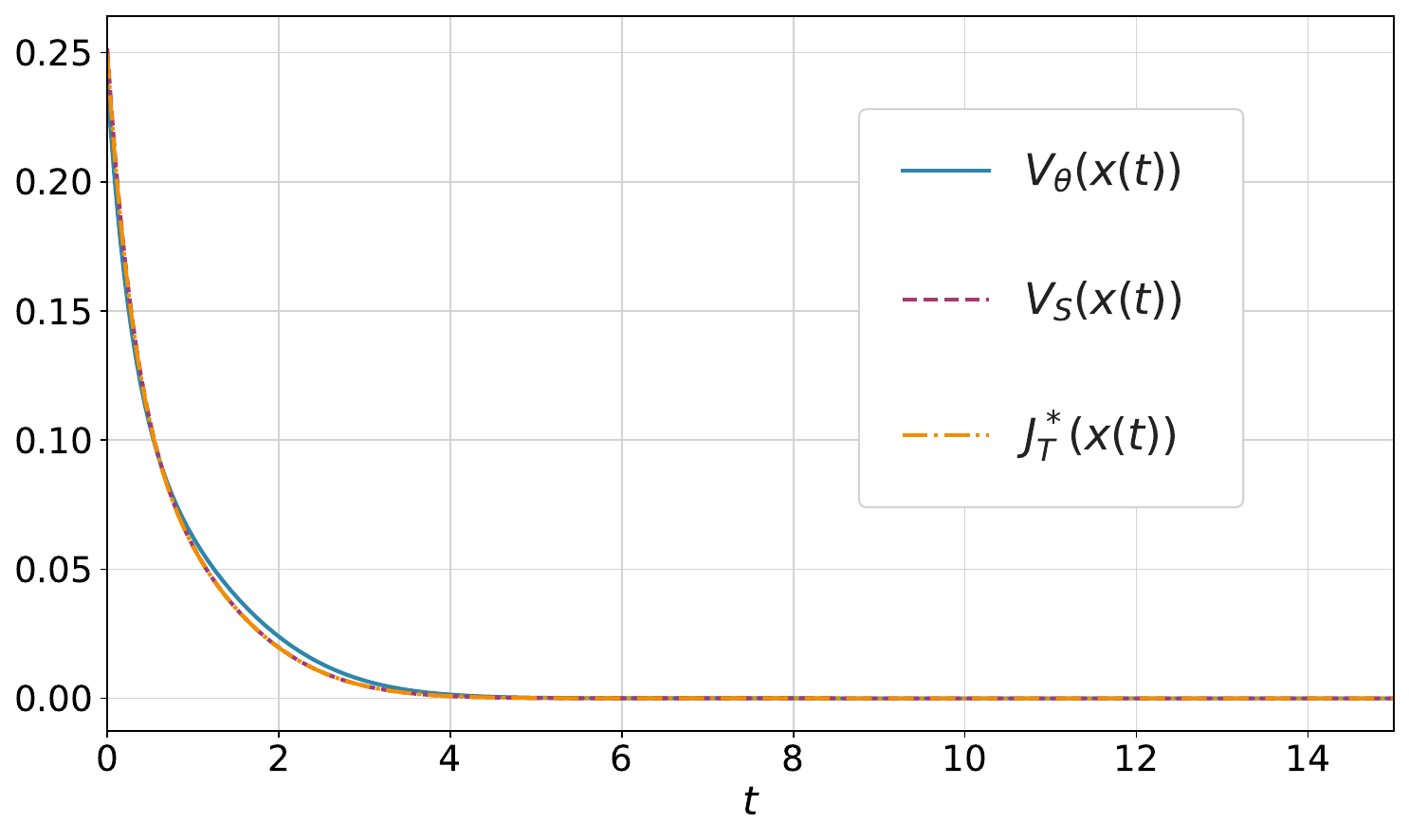}\qquad
\includegraphics[width=.38\linewidth]{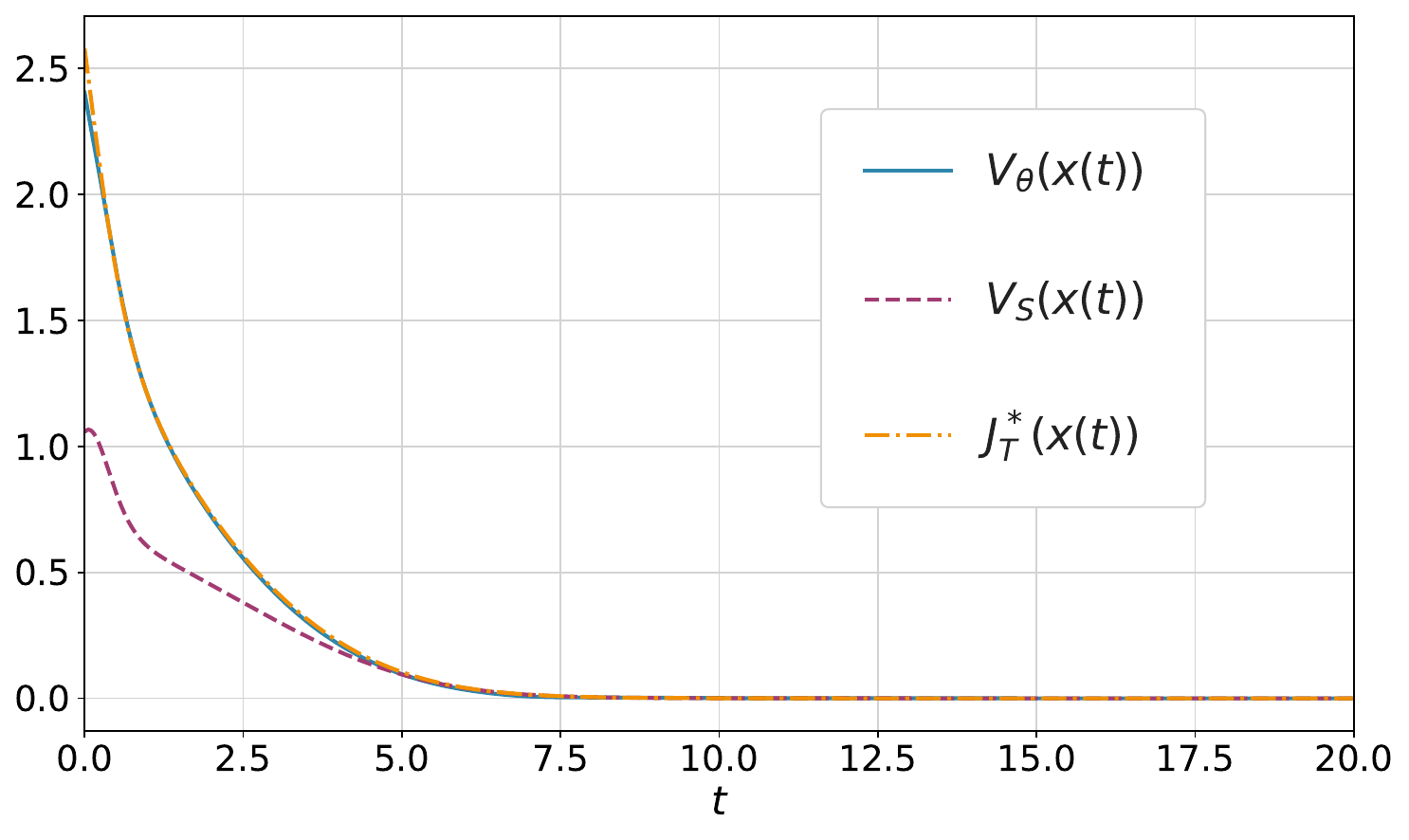}
\caption{
Comparison of the learned value function $V_\theta(x(t))$, the SDRE-based approximation $V_S(x(t))$, and the cost-to-go $J_T^*(x(t))$. Left: Cucker--Smale. Right: Allen--Cahn.
}
\label{fig:value}
\end{figure}

\noindent{\bf Result.}\;
We evaluate the trained NN controller and the baseline SDRE controller on test trajectories that were not used during training. We use equally spaced initial positions in $[0,0.2]$
and velocities in $[0.2,0.4]$ for Cucker--Smale ($T=15$), and $x_0(\xi)=1+\xi(1-\xi)$ for Allen--Cahn ($T=20$),
with a control-update interval of $0.01$ in both systems. Figures~\ref{fig:CuckerS_traj} and~\ref{fig:Allen-Cahn state} show that the NN controller stabilizes both systems to the origin, with a performance visually comparable to that of the SDRE controller.
The control signals in Figure~\ref{fig:control_signal} confirm this observation: for the Cucker--Smale model, the NN and SDRE controls are almost identical, and for the Allen--Cahn equation, the two controllers produce virtually overlapping signals.

Figure~\ref{fig:value} compares the value function learned by our method $V_\theta$, with the SDRE approximation $V_S$ and the cost-to-go $J_T^*$.
In the Cucker--Smale example, all three curves nearly coincide.
In the Allen--Cahn example, however, $V_S$ clearly deviates from $J_T^*$, while $V_\theta$ remains in close agreement with the cost-to-go.
This indicates that, at least for some nonlinear systems, the SDRE value function $V_S$ is not a faithful surrogate for the realized cost-to-go, and using it as a supervised label could introduce systematic bias.
By relying solely on control labels and the HJB residual, our method avoids this pitfall and recovers a more accurate value function.

Another advantage of the proposed method is its computational speed at inference.
Table~\ref{tab:time} reports the average wall-clock time required to evaluate a single control action.
In both examples, the NN controller is approximately 20-–30 times faster than the SDRE solver, which must solve a state-dependent Riccati equation at every step.
The speedup becomes even more significant in higher-dimensional settings or when the controller is called repeatedly inside a simulation loop.

\begin{table}[!hbt]\small
    \centering
    \caption{Average evaluation time (seconds) per state: SDRE controller vs.\ NN controller.}
    \label{tab:time}
    \begin{tabular}{l cc}
        \toprule
        Example & SDRE & $u_\theta$ \\
        \midrule
        Cucker--Smale ($d=40$) & 0.0063 & 0.0003 \\
        Allen--Cahn   ($d=51$) & 0.0115 & 0.0004 \\
        \bottomrule
    \end{tabular}
\end{table}

Additional experiments show that the HJB term improves
value--cost agreement, while control-only training can retain
comparable closed-loop costs. HJB-only training fails in the
tested cases. Loss and data-size ablations, loss-weight sensitivity,
label-noise robustness, OOD evaluation, tests under
perturbed dynamics, and full experimental details are
reported in Appendix~\ref{app:additional_experiments}.

\section{Conclusion}\label{sec:conclusion}



We presented a physics-informed method that learns optimal feedback laws from control labels only, without value function data. A neural network parameterizes the value function and is trained by jointly minimizing the HJB residual and the mismatch with SDRE-precomputed controls. Error bounds on the value function and its gradient are established through the viscous HJBE. Experiments on the Cucker--Smale and Allen--Cahn systems show that the resulting controller matches SDRE accuracy while being 20--30 times faster, and faithfully recovers the cost-to-go even when SDRE-based value approximations deviate. This method enables fast neural control with potential application in robotics and autonomous systems, though safety-critical deployment requires additional validation.

The proposed method has several limitations. First, it requires full knowledge of the system dynamics. Second, the method relies on precomputed SDRE controls, which can be computationally expensive in high dimensions. Third, the weighted loss formulation introduces a hyperparameter balancing different components, and optimal tuning may be problem-dependent. In addition, our theoretical guarantees are conditional.
The benchmarks fall outside the structural regime of the
policy-error theorem, and the finite trajectory data do not
certify the residual conditions required for global error guarantees.
The numerical results therefore provide empirical evidence within
the tested operating regions. Future work includes adaptive residual-based sampling and loss weighting, extending the policy-error analysis to underactuated systems and more general drift fields, potentially via projected-gradient and adjoint-weighted estimates, and learning with partially or fully unknown dynamics.

\section*{Acknowledgments}
This work was supported by grant of the National Natural Science
Foundation of China (No.~12501613). The authors declare no competing interests.

\bibliographystyle{elsarticle-num}
\bibliography{ref}

\appendix
\section{Technical Appendices and Supplementary Material}
We attached the corresponding proof details here.

\subsection{Viscosity solution method}
Since the weak regularity admitted by the HJBE, a type of weak solution called a viscosity solution was proposed.

\begin{definition}
V(x) is called a viscosity subsolution (supersolution) in $\Omega$ if for arbitrary test function $\phi\in C^1(\Omega)$, considering local maximum (minimum) $x_0$ of function $V - \phi$ satisfying $V(x_0) = \phi(x_0)$, the following inequality holds:
\[
\mathcal{H}[\phi](x_0)\le 0, \quad(\mathcal{H}[\phi](x_0)\ge 0)
\]
If a solution $V\in C^1(\Omega)$ is both a subsolution and a supersolution, we call it a viscosity solution.
\end{definition}

In the following, we shall introduce one type of technique to get a viscosity solution called the viscosity vanishing method. For the static HJBE \eqref{HJB_V}, it usually will have multiple solutions that satisfy the optimal condition. In order to solve this issue, the following viscosity form of HJB is considered
\begin{equation}
\lambda V^\varepsilon(x)
-\min_{u\in\mathbb{U}}
\left\{
\nabla V^\varepsilon(x)^\top(f(x)+B(x)u)+\ell(x,u)
\right\}
=\varepsilon\Delta V^\varepsilon(x).
\end{equation}

Substituting the closed-form optimal control \eqref{OC_u}, we get the expansion form.
\begin{equation} \label{viscous_HJB_eq}
\begin{aligned}
&\lambda V^\varepsilon(x) -x^\top Qx -  \nabla V^\varepsilon(x)^\top \cdot f(x) 
+ \frac{1}{4}\nabla V^\varepsilon(x)^\top B(x)R^{-1}B(x)^\top \nabla V^\varepsilon(x) = \varepsilon \Delta V^\varepsilon(x)
\end{aligned}
\end{equation} 
where $x\in\Omega$ and $\varepsilon\Delta V^\varepsilon$ denotes diffusion terms with small coefficient. Under some mild regularity condition about $f$ and $\ell$ \cite{evans2010adjoint}, the classical Bernstein method can be used to show \eqref{viscous_HJB_eq} has a unique classical solution $V^\varepsilon(x)$, which is a viscosity solution at the same time. A viscosity solution is easier to calculate by numerical methods due to the additional diffusion term. 

For the relation between solutions $V(x)$ and $V^\varepsilon(x)$, there is the following convergence result \cite{evans2010adjoint}.
\begin{equation}\label{Eq:V-Ve}
\| V(x) - V^\varepsilon(x)\|_{L^\infty(\R^n)}\le C\sqrt{\varepsilon}
\end{equation}
where $V$ denotes the viscosity solution of the original HJBE. 

For the case of a bounded state space $\Omega\neq\R^n$, the control policy should satisfy certain boundary conditions expressed by the inequality of the Hamiltonian.
\[
-\mathcal{H}(x, V, \nabla V+\alpha \hat{n})\le 0
\]
for any $\alpha\le 0$ and $x\in \partial\Omega$. Here $\hat{n}$ denotes the outer normal vector at the point $x$. For infinite-horizon HJBEs without state constraints, i.e., $\Omega = \R^n$, the “boundary condition” is typically replaced by a growth condition at infinity, because the domain is unbounded. Here’s a detailed explanation. The corresponding proof can be found in section \ref{pf: V growth}.

\subsection{Growth rate estimate}\label{pf: V growth}
\begin{proposition}\label{prop:V growth}
Consider the infinite-horizon discounted optimal control
problem (1). Assume that there exists an admissible control
$\bar{u}(\cdot)$ such that the corresponding trajectory
satisfies
\[
\|x(t;\bar{u})\|
\le C'(1+\|x\|)e^{Kt},
\qquad t\ge 0,
\]
for some constants $C'>0$ and $K>0$.
Furthermore, suppose that the running cost along this
trajectory satisfies
\[
0\le
\ell(x(t;\bar{u}),\bar{u}(t))
\le C\left(1+\|x(t;\bar{u})\|^p\right),
\]
for some constants $C>0$ and integer $p\ge 1$.
If $\lambda>Kp$, then the value function $V(x)$ satisfies
the polynomial growth estimate
\[
0\le V(x)\le C''(1+\|x\|^p),
\]
where $C''>0$ is independent of $x$.
\end{proposition}

\begin{proof}
By the definition of the value function and the
non-negativity of the running cost, we have
\begin{equation}
0\le V(x)
=\inf_{u(\cdot)\in\mathcal{U}}J(u(\cdot),x)
\le J(\bar{u}(\cdot),x),
\end{equation}
where $\bar{u}(\cdot)$ is the admissible control specified
in Proposition A.1.

Using the polynomial growth assumption on the running
cost and the exponential estimate for the corresponding
trajectory, we obtain
\begin{align}
\ell(x(t;\bar{u}),\bar{u}(t))
&\le C\left(1+\|x(t;\bar{u})\|^p\right)
\nonumber\\
&\le C\left(1+C'^p(1+\|x\|)^p e^{Kpt}\right)
\nonumber\\
&\le C_1\left(1+\|x\|^p e^{Kpt}
+e^{Kpt}\right),
\end{align}
where $C_1>0$ is independent of $x$ and $t$.

Consequently,
\begin{align}
V(x)
&\le J(\bar{u}(\cdot),x)
\nonumber\\
&=\int_0^\infty e^{-\lambda t}
\ell(x(t;\bar{u}),\bar{u}(t))\,dt
\nonumber\\
&\le C_1\int_0^\infty e^{-\lambda t}\,dt
\nonumber\\
&\quad+C_1(1+\|x\|^p)
\int_0^\infty e^{-(\lambda-Kp)t}\,dt.
\end{align}

Since $\lambda>Kp$, both integrals are finite, and hence
\[
V(x)
\le \frac{C_1}{\lambda}
+\frac{C_1}{\lambda-Kp}(1+\|x\|^p)
\le C''(1+\|x\|^p),
\]
for some constant $C''>0$ independent of $x$.
Therefore,
\[
0\le V(x)\le C''(1+\|x\|^p).
\]
This completes the proof.
\end{proof}

\subsection{Proof of Lemma \ref{Lem:SDPP}}\label{pf:Lem:SDPP}
After obsorbing the extra residual $r_\vare$ into the minimization, we can regard the NN value approximation $V_\theta(x)$ as solution of modified viscosity HJB with modified running cost $\tilde \ell_\vare$. It corresponds to a stochastic optimal control optimization in which NN $V_\theta$ is exactly the optimal value of a stochastic process. Based on the classical stochastic optimal control theory \cite{fleming2012deterministic}, the Stochastic Dynamical Programming Principle (SDPP) will be obtained.

\hfill$\blacksquare$

\subsection{Proof of Theorem \ref{prop:VS-Ve}}\label{pf:prop:VS-Ve}

For every $T\ge 0$, we write down the SDPP formulation for the true viscosity solution.
\begin{equation}
V^\varepsilon(x)
=
\inf_{u\in\mathcal U}
\mathbb E_x\!\left[
\int_0^T e^{-\lambda t}\ell(X_t,u_t)\,dt
+
e^{-\lambda T}V^\varepsilon(X_T)
\right].
\end{equation}
In particular, for the optimal control $u^*$,
\begin{equation}\label{Eq:Ve_DPP}
V^\varepsilon(x)
=
\mathbb E_x\!\left[
\int_0^T e^{-\lambda t}\ell(X_t^*,u_t^*)\,dt
+
e^{-\lambda T}V^\varepsilon(X_T^*)
\right].
\end{equation}
On the other hand, since $V_\theta$ solves the modified viscosity form HJB with running cost $\tilde\ell_\varepsilon$, its SDPP yields the following form
\[
V_\theta(x)
=
\inf_{u\in\mathcal U}
\mathbb E_x\!\left[
\int_0^T e^{-\lambda t}\tilde\ell_\varepsilon(X_t,u_t)\,dt
+
e^{-\lambda T}V_\theta(X_T)
\right].
\]
Evaluating this identity along the control $u^*$, we get
\begin{equation}\label{Eq:VS-1}
    V_\theta(x)
\le
\mathbb E_x\!\left[
\int_0^T e^{-\lambda t}\tilde\ell_\varepsilon(X_t^*,u_t^*)\,dt
+
e^{-\lambda T}V_\theta(X_T^*)
\right].
\end{equation}

Using the relation between running cost $\tilde\ell_\varepsilon(x,u)=\ell(x,u)+r_\varepsilon(x)$. Inequality \eqref{Eq:VS-1}  becomes
\begin{equation}\label{Eq:VS-2}
    V_\theta(x)
\le
\mathbb E_x\!\left[
\int_0^T e^{-\lambda t}\bigl(\ell(X_t^*,u_t^*)+r_\varepsilon(X_t^*)\bigr)\,dt
+
e^{-\lambda T}V_\theta(X_T^*)
\right].
\end{equation}

Subtracting (\ref{Eq:Ve_DPP}) from (\ref{Eq:VS-2}), we obtain
\[
V_\theta(x)-V^\varepsilon(x)
\le
\mathbb E_x\!\left[
\int_0^T e^{-\lambda t}r_\varepsilon(X_t^*)\,dt
\right]
+
\mathbb E_x\!\left[
e^{-\lambda T}\bigl(V_\theta(X_T^*)-V^\varepsilon(X_T^*)\bigr)
\right].
\]
i.e.
\[
V_\theta(x)-V^\varepsilon(x)
\le
\mathbb E_x\!\left[
\int_0^T e^{-\lambda t}|r_\varepsilon(X_t^*)|\,dt
\right]
+
\mathbb E_x\!\left[
e^{-\lambda T}\bigl|V_\theta(X_T^*)-V^\varepsilon(X_T^*)\bigr|
\right].
\]

Passing to the limit as $T\to\infty$, and using the terminal condition assumption, we arrive at
\begin{equation}\label{Eq:VS-Ve}
V_\theta(x)-V^\varepsilon(x)
\le
\mathbb E_x\!\int_0^\infty e^{-\lambda t}|r_\varepsilon(X_t^*)|\,dt.    
\end{equation}

To obtain the reverse estimate, we apply the SDPP for $V_\theta$ along its own optimal control $\widetilde u$:
\begin{equation}\label{Eq:VS-3}
    V_\theta(x)
=
\mathbb E_x\!\left[
\int_0^T e^{-\lambda t}\tilde\ell_\varepsilon(\widetilde X_t,\widetilde u_t)\,dt
+
e^{-\lambda T}V_\theta(\widetilde X_T)
\right].
\end{equation}

Since $V^\varepsilon$ is the infimum of the original cost, evaluating its SDPP along the same control $\widetilde u$ gives
\begin{equation}\label{Eq:Ve-1}
V^\varepsilon(x)
\le
\mathbb E_x\!\left[
\int_0^T e^{-\lambda t}\ell(\widetilde X_t,\widetilde u_t)\,dt
+
e^{-\lambda T}V^\varepsilon(\widetilde X_T)
\right].    
\end{equation}

Subtracting (\ref{Eq:VS-3}) from (\ref{Eq:Ve-1}), and using again $\tilde\ell_\varepsilon=\ell+r_\varepsilon$, we obtain
\[
V^\varepsilon(x)-V_\theta(x)
\le
\mathbb E_x\!\left[
\int_0^T e^{-\lambda t}|r_\varepsilon(\widetilde X_t)|\,dt
\right]
+
\mathbb E_x\!\left[
e^{-\lambda T}\bigl|V^\varepsilon(\widetilde X_T)-V_\theta(\widetilde X_T)\bigr|
\right].
\]
Letting $T\to\infty$, we deduce
\begin{equation}\label{Eq:Ve-VS}
V^\varepsilon(x)-V_\theta(x)
\le
\mathbb E_x\!\int_0^\infty e^{-\lambda t}|r_\varepsilon(\widetilde X_t)|\,dt.    
\end{equation}
Combining (\ref{Eq:Ve-VS}) and (\ref{Eq:VS-Ve}), we conclude that
\[
|V_\theta(x)-V^\varepsilon(x)|
\le
\max\left\{
\mathbb E_x\!\int_0^\infty e^{-\lambda t}|r_\varepsilon(X_t^*)|\,dt,\,
\mathbb E_x\!\int_0^\infty e^{-\lambda t}|r_\varepsilon(\widetilde X_t)|\,dt
\right\}.
\]
Finally, if $\lambda>0$ and $r_\varepsilon\in L^\infty$, then for any admissible state process $X_t$,
\[
\mathbb E_x\!\int_0^\infty e^{-\lambda t}|r_\varepsilon(X_t)|\,dt
\le
\|r_\varepsilon\|_{L^\infty}\int_0^\infty e^{-\lambda t}\,dt
=
\frac{1}{\lambda}\|r_\varepsilon\|_{L^\infty}.
\]
Therefore,
\[
|V_\theta(x)-V^\varepsilon(x)|
\le
\frac{1}{\lambda}\|r_\varepsilon\|_{L^\infty}.
\]
This completes the proof.\hfill$\blacksquare$

\subsection{Liouville-type inequality}

\begin{lemma}\label{Lem-energy-estimate}
Let $B_r:=\{x\in \mathbb{R}^n: |x|\le r\}$ be a bounded ball. Assume that the control matrix $B(x)$ has full rank and  $B(x)R^{-1}B(x)^\top=I$. Assume the drift function is of the gradient form $f=\nabla \phi$. Denote $\Phi(x):=\frac12 V^\varepsilon(x)-\phi(x)$ and $\rho(x):=\exp\left(-\frac{\Phi(x)}{\varepsilon}\right)$. Let $R_\varepsilon(x):=r_\varepsilon[V_\theta](x)+\lambda w(x)$ be the error source term derived from the viscosity form HJB and decay term. Then we have the following energy estimate for weighted energy.
\begin{equation}\label{Eq:E(r)}
\frac14 E(r)
\le
\varepsilon \sqrt{E'(r)M'(r)}
+
H(r).
\end{equation}
where $E(r):=\int_{B_r} \rho |\nabla w|^2, M(r):=\int_{B_r} \rho, H(r):=\int_{B_r} \rho |R_\varepsilon|$.
\end{lemma}

\begin{proof}
Subtract \eqref{viscous_HJB_eq} and \eqref{Eq:Vtheta}, we re-express the equation in terms of the error.
\begin{equation}\label{Eq:Vtheta-Ve}
\lambda w-\nabla w^{\top} f
+\frac12 \nabla V^\varepsilon \cdot \nabla w
-\frac14 |\nabla w|^2
=
\varepsilon \Delta w-r_\varepsilon[V_\theta].
\end{equation}

For simplicity, we assume that the control matrix $B(x)$ has full rank and that $B(x)R^{-1}B(x)^\top=I$. We further assume that $f=\nabla \phi$. Define
\[
\Phi(x):=\frac12 V^\varepsilon(x)-\phi(x),
\qquad
R_\varepsilon(x):=r_\varepsilon[V_\theta](x)+\lambda w(x).
\]
Then \eqref{Eq:Vtheta-Ve} can be rewritten as
\[
-\varepsilon \Delta w+\nabla \Phi \cdot \nabla w
-\frac14 |\nabla w|^2
=
-R_\varepsilon .
\]

Let
\begin{equation}\label{Eq:rho}
\rho(x):=\exp\left(-\frac{\Phi(x)}{\varepsilon}\right).
\end{equation}
Since
\[
\rho\left(-\varepsilon \Delta w+\nabla \Phi \cdot \nabla w\right)
=
-\varepsilon \nabla \cdot (\rho \nabla w),
\]
multiplying \eqref{Eq:rho} by $\rho$ gives
\[
-\varepsilon \nabla \cdot (\rho \nabla w)
-\frac14 \rho |\nabla w|^2
=
-\rho R_\varepsilon .
\]

Integrating over $B_r:=\{x\in \mathbb{R}^n: |x|\le r\}$, 
and applying the divergence theorem, we obtain
\[
-\varepsilon \int_{\partial B_r} \rho \nabla w \cdot \nu
-\frac14 \int_{B_r} \rho |\nabla w|^2
=
-\int_{B_r} \rho R_\varepsilon ,
\]
where $\nu$ denotes the unit outward normal vector on $\partial B_r$.
Hence,
\[
\frac14 \int_{B_r} \rho |\nabla w|^2
=
\int_{B_r} \rho R_\varepsilon
-
\varepsilon \int_{\partial B_r} \rho \nabla w \cdot \nu .
\]

Taking absolute values and applying the Cauchy-Schwarz inequality to the boundary term yields
\begin{equation}\label{Eq:grad_w}
\frac14 \int_{B_r} \rho |\nabla w|^2
\le
\varepsilon
\left(\int_{\partial B_r} \rho |\nabla w|^2\right)^{1/2}
\left(\int_{\partial B_r} \rho\right)^{1/2}
+
\int_{B_r} \rho |R_\varepsilon|.
\end{equation}

Now define
\[
E(r):=\int_{B_r} \rho |\nabla w|^2,
\qquad
M(r):=\int_{B_r} \rho,
\qquad
H(r):=\int_{B_r} \rho |R_\varepsilon|.
\]
Then, for almost every $r>0$,
\[
E'(r)=\int_{\partial B_r} \rho |\nabla w|^2,
\qquad
M'(r)=\int_{\partial B_r} \rho .
\]

Therefore, estimate \eqref{Eq:grad_w} becomes the following energy estimate
\[
\frac14 E(r)
\le
\varepsilon \sqrt{E'(r)M'(r)}
+
H(r).
\]
Equivalently,
\begin{equation}
\left(\frac14 E(r)-H(r)\right)_+^2
\le
\varepsilon^2 E'(r)M'(r),
\end{equation}
where $(a)_+:=\max\{a,0\}$.
\end{proof}

\subsection{Proof of Theorem \ref{Thm-error-gradient-estimate}}\label{pf:Thm-error-gradient-estimate}

Recall we have derived error energy estimate (\ref{Eq:E(r)}):
\begin{equation}
\frac14 E(r)
\le
\varepsilon \sqrt{E'(r)M'(r)}
+
H(r).
\label{eq:key-estimate}
\end{equation}

By the statement assumptions for $M(r)$ and $H(r)$, we have
\[
M'(r)\le C,
\qquad
H(r)\le H_\infty\le \frac{\delta}{8}.
\]
Thus \eqref{eq:key-estimate} implies
\begin{equation}
\frac14 E(r)
\le
\varepsilon \sqrt{C E'(r)}
+
H_\infty .
\label{eq:key-estimate-C}
\end{equation}

We now prove the desired estimate by contradiction.
Suppose that there exists $r_0>0$ such that
\begin{equation}
E(r_0)>\delta .
\end{equation}
Since $E(r)$ is nondecreasing in $r$, it follows that
\begin{equation}
E(r)\ge E(r_0)>\delta
\qquad
\text{for all } r\ge r_0.
\end{equation}
Hence,
\begin{equation}
H_\infty
\le
\frac{\delta}{8}
<
\frac{E(r)}{8}
\qquad
\text{for all } r\ge r_0.
\end{equation}
Consequently,
\begin{equation}
\frac14 E(r)-H_\infty
\ge
\frac18 E(r).
\end{equation}
From \eqref{eq:key-estimate-C}, we then obtain
\begin{equation}
\frac18 E(r)
\le
\varepsilon \sqrt{C E'(r)}
\qquad
\text{for all } r\ge r_0.
\end{equation}
Squaring both sides gives
\begin{equation}
E(r)^2
\le
64\varepsilon^2 C E'(r).
\end{equation}
Equivalently,
\begin{equation}
\frac{1}{64\varepsilon^2 C}
\le
\frac{E'(r)}{E(r)^2}.
\end{equation}

Integrating this inequality over $(r_0,R)$, we obtain
\begin{equation}
\frac{R-r_0}{64\varepsilon^2 C}
\le
\int_{r_0}^{R}
\frac{E'(r)}{E(r)^2}\,dr .
\end{equation}
Since
\begin{equation}
\int_{r_0}^{R}
\frac{E'(r)}{E(r)^2}\,dr
=
\frac{1}{E(r_0)}-\frac{1}{E(R)},
\end{equation}
we have
\begin{equation}
\frac{R-r_0}{64\varepsilon^2 C}
\le
\frac{1}{E(r_0)}-\frac{1}{E(R)}
\le
\frac{1}{E(r_0)}.
\end{equation}
Letting $R\to\infty$, the left-hand side diverges, while the right-hand side remains finite.
This contradiction shows that no such $r_0$ can exist.

Therefore,
\begin{equation}
E(r)\le \delta
\qquad
\text{for every } r>0.
\end{equation}
That is,
\begin{equation}
\int_{B_r}
\rho |\nabla w|^2
\le
\delta.
\end{equation}
Substituting back $\rho=e^{-\Phi/\varepsilon}$, we conclude that
\begin{equation}
\int_{B_r}
\exp\left(-\frac{\Phi}{\varepsilon}\right)
|\nabla w|^2
\le
\delta,
\end{equation}
which proves the desired weighted gradient estimate. By applying Young's inequality, we have
\[
\int_{B_r}\rho |\nabla V_\theta - \nabla V|^2\le 2 \int_{B_r}\rho |\nabla V_\theta - \nabla V^\vare|^2 + 2 \int_{B_r}\rho |\nabla V^\vare - \nabla V|^2
\]
It follows that
\begin{equation}
\int_{B_r}\rho|\nabla V_\theta - \nabla V|^2 \le 2\delta + 2\eta_\vare(r)
\end{equation}

In the following, we shall focus on the approximation of the control function by NN. The neural-network-induced feedback control is
\begin{equation}
u_\theta(x)
=
-\frac12 R^{-1}B(x)^{\top}\nabla V_\theta(x).
\end{equation}
Similarly, the viscosity optimal control associated with $V^\varepsilon$ is
\begin{equation}
u^\varepsilon(x)
=
-\frac12 R^{-1}B(x)^{\top}\nabla V^\varepsilon(x).
\end{equation}
and true control function 
\[
u(x)=-\frac{1}{2}R^{-1}B(x)^\top\nabla V(x)
\]

Therefore,
\begin{align}
u_\theta(x)-u(x)
&=
-\frac12 R^{-1}B(x)^{\top}
\left(
\nabla V_\theta(x)-\nabla V(x)
\right).
\end{align}

Taking norms gives the point-wise estimate
\begin{equation}
|u_\theta(x)-u(x)|
\le
\frac12
\left\|
R^{-1}B(x)^{\top}
\right\|
|\nabla V_\theta - \nabla V|.
\end{equation}

Assume that
\begin{equation}
\left\|
R^{-1}B(x)^{\top}
\right\|
\le K_r,
\qquad x\in B_r .
\end{equation}
Then
\begin{equation}
\begin{aligned}
    \int_{B_r}\rho |u_\theta - u|^2\le& \frac{K_r^2}{4}\int_{B_r}\rho |\nabla V_\theta - \nabla V|^2\\
    \le& \frac{K_r^2}{2}(\delta + \eta_\vare)
\end{aligned}
\end{equation}
This gives the corresponding optimal control error estimate.
\hfill$\blacksquare$

\subsection{Proof of Theorem \ref{Thm-policy-estimate}}\label{pf:Thm-policy-estimate}

Recall that $V^\varepsilon$ satisfies the viscous HJB equation
\begin{equation}
\label{eq:viscous_hjb_policy}
\lambda V^\varepsilon
-x^\top Qx
-f^\top\nabla V^\varepsilon
+\frac14(\nabla V^\varepsilon)^\top
G\nabla V^\varepsilon
-\varepsilon\Delta V^\varepsilon=0,
\end{equation}
where $G(x)=B(x)R^{-1}B(x)^\top$. Similarly, by the definition of the viscous HJB
residual, $V_\theta$ satisfies
\begin{equation}
\label{eq:nn_residual_policy}
\lambda V_\theta
-x^\top Qx
-f^\top\nabla V_\theta
+\frac14(\nabla V_\theta)^\top
G\nabla V_\theta
-\varepsilon\Delta V_\theta
=r_\varepsilon[V_\theta].
\end{equation}

Let $w=V^\varepsilon-V_\theta$ then subtracting \eqref{eq:nn_residual_policy} from
\eqref{eq:viscous_hjb_policy} gives
\begin{align}
\lambda w-f^\top\nabla w
&+\frac14\Big(
(\nabla V^\varepsilon)^\top G\nabla V^\varepsilon
-(\nabla V_\theta)^\top G\nabla V_\theta
\Big)
\nonumber\\
&-\varepsilon\Delta w
=-r_\varepsilon[V_\theta].
\end{align}

Since $\nabla V_\theta=\nabla V^\varepsilon-\nabla w$,
the quadratic difference can be written as
\[
(\nabla V^\varepsilon)^\top G\nabla V^\varepsilon
-(\nabla V_\theta)^\top G\nabla V_\theta =
2(\nabla V^\varepsilon)^\top G\nabla w
-(\nabla w)^\top G\nabla w.
\]
Consequently, error $w$ will satisfy the following:
\begin{equation}
\label{eq:nonlinear_error_policy}
-\varepsilon\Delta w
+b_0^\top\nabla w
-\frac14(\nabla w)^\top G\nabla w
+\lambda w
=-r_\varepsilon[V_\theta],
\end{equation}
where $b_0:=-f+\frac12G\nabla V^\varepsilon$. Next, observe that $u_\theta-u^\varepsilon
=\frac12R^{-1}B^\top\nabla w.$ Therefore, it is directly calculated that $\frac14(\nabla w)^\top G\nabla w = |u_\theta-u^\varepsilon|_R^2$ where $|v|_R^2=v^\top Rv$.

Let $q(x)=|u_\theta-u^\varepsilon|_R^2\geq0$ and \eqref{eq:nonlinear_error_policy} then becomes
\begin{equation}
\label{eq:error_q_identity}
q=r_\varepsilon[V_\theta]
-\varepsilon\Delta w
+b_0^\top\nabla w+\lambda w.
\end{equation}

For any $B_r\Subset B_{r_1}\Subset\Omega$,
choose a smooth cutoff function
$\chi\in C_c^\infty(B_{r_1})$ satisfying
\[
0\leq\chi\leq1,
\qquad
\chi=1 \quad\text{on }B_r.
\]
Multiplying \eqref{eq:error_q_identity} on both sides by $\chi$
and integrating over $B_{r_1}$, we then use the integration by parts to get


\begin{equation}
\label{eq:cutoff_policy_identity}
\int_{B_{r_1}}\chi q\,dx
=
\int_{B_{r_1}}\chi r_\varepsilon[V_\theta]\,dx+\int_{B_{r_1}}w
\left(
\lambda\chi-\varepsilon\Delta\chi
-\nabla\cdot(\chi b_0)
\right)dx.
\end{equation}

Since $q\geq0$ and $\chi=1$ on $B_r$, it follows that $\int_{B_r}q\,dx
\leq
\int_{B_{r_1}}\chi q\,dx$. Using \eqref{eq:cutoff_policy_identity},
we conclude that
\begin{equation}
\label{eq:weighted_control_intermediate}
\begin{aligned}
\int_{B_r}|u_\theta-u^\varepsilon|_R^2\,dx
\leq{}&
\|r_\varepsilon[V_\theta]\|_{L^1(B_{r_1})}+C_{\varepsilon,r,r_1}
\|V^\varepsilon-V_\theta\|_{L^1(B_{r_1})}.
\end{aligned}
\end{equation}
where $C_{\varepsilon,r,r_1}
=
\left\|
\lambda\chi-\varepsilon\Delta\chi
-\nabla\cdot(\chi b_0)
\right\|_{L^\infty(B_{r_1})}$ which is independent of $V_\theta$.

Since $R\succ0$, we have
\[
|v|^2
\leq\frac{1}{\lambda_{\min}(R)}|v|_R^2.
\]
Thus,
\[
\begin{aligned}
\|u_\theta-u^\varepsilon\|_{L^2(B_r)}^2
\leq
\frac{1}{\lambda_{\min}(R)}
\Big(
&\|r_\varepsilon[V_\theta]\|_{L^1(B_{r_1})}\\
&+C_{\varepsilon,r,r_1}
\|V^\varepsilon-V_\theta\|_{L^1(B_{r_1})}
\Big).
\end{aligned}
\]

Finally, by the triangle inequality,
\[
|u_\theta-u|^2
\leq
2|u_\theta-u^\varepsilon|^2
+2|u^\varepsilon-u|^2.
\]
Integrating over $B_r$ and using
\[
\|u^\varepsilon-u\|_{L^2(B_r)}^2
\leq\eta_\varepsilon(r)
\]
gives
\[
\begin{aligned}
\|u_\theta-u\|_{L^2(B_r)}^2
\leq
\frac{2}{\lambda_{\min}(R)}
\Big(
&\|r_\varepsilon[V_\theta]\|_{L^1(B_{r_1})}+C_{\varepsilon,r,r_1}
\|V^\varepsilon-V_\theta\|_{L^1(B_{r_1})}
\Big)
+2\eta_\varepsilon(r).
\end{aligned}
\]
This completes the proof.

\hfill$\blacksquare$

\section{Additional Experimental Details and Results}
\label{app:additional_experiments}

\subsection{Ablations and Loss-Weight Sensitivity}
\label{app:ablation_gamma}
\label{app:experimental_protocol}
Table~\ref{tab:baseline_comparison} compares the proposed
method with two baselines: the same value-function network
trained using $\mathcal L_c$ alone, with controls still
derived from its gradient, and a feedforward policy network
that we denote by \emph{PolicyNet}.
PolicyNet directly maps states to controls and minimizes
the mean squared error between its output and the SDRE labels; it does not estimate
a value function. For this comparison, all three neural methods follow
the hidden-layer configurations and optimization settings
of the main-text experiments in Section~\ref{sec:experiment}.
They share a common 80/20 pointwise training/held-out
partition within each benchmark: 16,000/4,000 samples
for Cucker--Smale and 8,000/2,000 for Allen--Cahn. 
Thus, for Allen--Cahn, this 80/20 split differs from the 90/10
split used in the main-text experiments.

Table~\ref{tab:ablation_gamma} also examines sensitivity to
$\gamma$ in the objective
$\mathcal L=\mathcal L_{\rm HJB}+\gamma\mathcal L_c$.
Within each system, we vary $\gamma$ while fixing the
dataset, 80/20 pointwise training/held-out partition,
network architecture, initialization seed, and nominal
test initial condition. Here, $\gamma=0$ denotes \emph{HJB-only} training,
while $\gamma=\infty$ denotes \emph{$u$-only} training
using $\mathcal L_c$ alone. All comparisons use the nominal test initial conditions,
horizons, and control-update interval specified in
Section~\ref{sec:experiment}. 

We report the closed-loop cost ratio
$C=J_{\rm NN}/J_{\rm SDRE}$, with both costs evaluated
from the same initial condition over the same horizon.
To quantify value--cost agreement, all value networks
are evaluated along the same SDRE reference trajectory
$x_S(t)$. The SDRE controller's realized remaining cost is
\[
\widehat J_S(t)
=\int_t^T \ell\bigl(x_S(\tau),u_S(x_S(\tau))\bigr)\,d\tau.
\]
We define the normalized area discrepancy between the
predicted value and this reference cost-to-go as
\begin{equation}
E_{\rm AUC}
=
\frac{\displaystyle\int_0^T
\left|V_\theta(x_S(t))-\widehat J_S(t)\right|\,dt}
{\displaystyle\int_0^T|\widehat J_S(t)|\,dt}.
\label{eq:eauc}
\end{equation}
This dimensionless metric measures the integrated absolute
discrepancy relative to the area under the reference
cost-to-go curve; lower values indicate closer agreement.
It assesses agreement with the finite-horizon SDRE
cost-to-go along the reference trajectory.

\begin{table}[tbp]
\centering
\setlength{\tabcolsep}{4pt}
\caption{Ablations and loss-weight sensitivity.
Both studies use an 80/20 training/held-out split
for each system.
Each entry uses one trained model and the nominal test
initial condition for that system.
$C=J/J_{\rm SDRE}$, and $E_{\rm AUC}$ is evaluated
on the SDRE reference trajectory defined in \eqref{eq:eauc}.
HJB-only results are retained as failure diagnostics.
Dashes indicate that no value metric is reported.}
\label{tab:ablation_gamma}
\label{tab:baseline_comparison}
\begin{tabular}{lrrrr}
\toprule
& \multicolumn{2}{c}{Cucker--Smale}
& \multicolumn{2}{c}{Allen--Cahn} \\
Setting
& $C$ & $E_{\rm AUC}$
& $C$ & $E_{\rm AUC}$ \\
\midrule
$\gamma=0$ (HJB-only)
& 4710.4202 & 0.8676 & 13.5656 & 0.6801 \\
$\gamma=10^{-3}$
& 1.1801 & 0.4705 & 1.0133 & 0.0113 \\
$\gamma=10^{-2}$
& 1.0186 & 0.0610 & 1.0075 & 0.0137 \\
$\gamma=10^{-1}$
& 1.0161 & 0.0343 & 1.0041 & 0.0094 \\
$\gamma=1$ (Ours)
& 1.0159 & 0.0449 & 1.0031 & 0.0361 \\
$\gamma=10$
& 1.0053 & 0.0671 & 1.0012 & 0.0738 \\
$\gamma=10^{2}$
& 1.0102 & 0.0382 & 1.0020 & 0.2099 \\
$\gamma=10^{3}$
& 1.0119 & 0.3563 & 0.9993 & 0.3950 \\
$\gamma=\infty$ ($u$-only)
& 1.0077 & 0.5251 & 0.9987 & 0.9930 \\
\midrule
PolicyNet
& 1.0578 & -- & 0.9983 & -- \\
SDRE reference
& 1.0000 & -- & 1.0000 & -- \\
\bottomrule
\end{tabular}
\end{table}

The three control-supervised methods achieve comparable nominal closed-loop costs. The proposed method additionally yields a value estimate that closely tracks the realized SDRE cost-to-go, whereas the $u$-only value network exhibits substantially larger value--cost discrepancies. The $\gamma$ study further shows that large control-loss weights can preserve near-reference control costs despite large value discrepancies, highlighting the role of the HJB term in improving value--cost consistency. HJB-only training yields substantially higher closed-loop costs in both systems, underscoring the importance of control supervision.

\subsection{Repeated Label-Noise Experiments}
\label{app:noise}
We add independent zero-mean Gaussian noise to each
SDRE control-label entry, with a standard deviation
equal to $r$ times the sample standard deviation of
all original control-label entries, where
$r\in\{0.10,0.25,0.50,0.75\}$. Five noise realizations use seeds 2026--2030, while
initialization and split seeds remain fixed at 42.
Within each system, all methods share the states,
noisy labels, 80/20 pointwise training/held-out partition,
and test initial conditions.
We set $\gamma=1$ for the proposed method.

We use method-specific network configurations.
For Cucker--Smale, the proposed method and $u$-only use two hidden
layers of 128 and 256 units, respectively; PolicyNet
uses three hidden layers of 512 units.
For Allen--Cahn, both value networks use two hidden layers of
128 units, and PolicyNet uses three of 512 units.
All networks use $\tanh$ activations and Adam with
batch size 4096.
For Cucker--Smale, we use 4000 epochs and an initial learning rate of
$8\times10^{-4}$; for Allen--Cahn, we use 3000 epochs and
$2\times10^{-3}$.

Each trained model is evaluated on five common test
initial conditions: uniform samples from $[-1,1]^{40}$
for Cucker--Smale and Fourier profiles from the
family in Section~\ref{sec:experiment} for Allen--Cahn. For each noise realization, we average
$C=J_{\rm NN}/J_{\rm SDRE}$ over the five test conditions,
then report the mean and sample standard deviation
across the five realization averages. These results in Table~\ref{tab:noise_robustness} suggest a problem-dependent benefit under
noisy supervision, with more pronounced improvements
at higher noise levels, particularly for Allen--Cahn.
The observations are specific to the tested systems
and network configurations and do not establish a
uniform advantage across problems or noise conditions.

\begin{table}[tbp]
\centering
\setlength{\tabcolsep}{4pt}
\caption{Label-noise robustness under the reported
method-specific network configurations.
Entries report normalized closed-loop cost,
mean $\pm$ sample standard deviation across five
noise-realization averages, each based on five common
test initial conditions.}
\label{tab:noise_robustness}
\begin{tabular}{lcccc}
\toprule
System & $r$ & $u+\mathrm{HJB}$ & $u$-only & PolicyNet \\
\midrule
Cucker--Smale & 0.10
& $1.0141\pm0.0019$ & $1.0494\pm0.0104$ & $1.0276\pm0.0018$ \\
Cucker--Smale & 0.25
& $1.0190\pm0.0032$ & $1.0586\pm0.0142$ & $1.0402\pm0.0032$ \\
Cucker--Smale & 0.50
& $1.0301\pm0.0059$ & $1.0843\pm0.0222$ & $1.0818\pm0.0068$ \\
Cucker--Smale & 0.75
& $1.0447\pm0.0122$ & $1.1181\pm0.0156$ & $1.1649\pm0.0097$ \\
\midrule
Allen--Cahn & 0.10
& $1.0030\pm0.0020$ & $1.0079\pm0.0043$ & $1.0055\pm0.0032$ \\
Allen--Cahn & 0.25
& $1.0128\pm0.0035$ & $1.0206\pm0.0156$ & $1.0670\pm0.0587$ \\
Allen--Cahn & 0.50
& $1.0635\pm0.0210$ & $1.1313\pm0.1535$ & $1.6970\pm0.4303$ \\
Allen--Cahn & 0.75
& $1.2417\pm0.1151$ & $1.3815\pm0.3132$ & $3.9040\pm1.3282$ \\
\bottomrule
\end{tabular}
\end{table}

\subsection{OOD and Long-Horizon Evaluation}
\label{app:ood}

We evaluate the original main-text $\gamma=1$ models
without retraining or modifying their feedback laws.
For Allen--Cahn, this is the model trained with the original
90/10 split. We extend the nominal horizons of 15 for
Cucker--Smale and 20 for Allen--Cahn to $T=50$, and test five initial
conditions per system at each scale
$s\in\{1,1.25,1.5,2\}$.

For Cucker--Smale, we independently sample five pairs of position
and velocity vectors $y,v\in[-1,1]^{20}$.
Each block is normalized separately, and the initial
state is constructed as
\[
x_0
=s\begin{bmatrix}
y/\|y\|_\infty\\
v/\|v\|_\infty
\end{bmatrix}.
\]
Here, $\|\cdot\|_\infty$ is the largest absolute component.
Thus, the normalized position and velocity blocks each
have maximum absolute entry one before scaling.
The same five normalized directions are reused across
all scales. Since the training initial conditions lie
in $[-1,1]^{40}$, the cases with $s=1$ lie on its boundary,
whereas those with $s>1$ lie outside it.

For Allen--Cahn, the test profiles are
\[
x_0(\xi)
=s\bigl[1+\xi(1-\xi)+\delta(\xi)\bigr].
\]
One profile uses $\delta=0$, while the other four use
smooth Fourier perturbations. Each perturbation is
centered to have zero mean over the 51 spatial grid
points before being added to the baseline profile.
The same five profiles are reused across scales,
with the entire profile multiplied by $s$. At all tested scales, the test-profile means lie outside
the empirical range of the training initial-profile means.

A rollout is classified as successful if it reaches
$T=50$ and its sampled states satisfy
\[
\|x(t)\|_2
\le \max\bigl\{10^{-3},\,0.05\|x_0\|_2\bigr\}
\qquad
\text{for all sampled }t\in[40,50].
\]
This requires the trajectory to enter a neighborhood
of the origin by time 40 and remain there until time 50. 


\begin{table}[tbp]
\centering
\setlength{\tabcolsep}{4pt}
\caption{OOD evaluation over a horizon of $T=50$.
Cost ratios relative to SDRE are medians over successful
cases only. Final state norms and absolute HJB residuals
are averaged over all five cases, including failures.}
\label{tab:ood_long_horizon}
\begin{tabular}{lccrrr}
\toprule
System & Scale & Successes
& \shortstack{Median\\cost ratio}
& \shortstack{Mean final\\state norm}
& \shortstack{Mean absolute\\HJB residual} \\
\midrule
Allen--Cahn & 1.00 & 5/5 & 1.0022 & 0.0042535 & 0.004296 \\
Allen--Cahn & 1.25 & 5/5 & 1.0046 & 0.0042535 & 0.011913 \\
Allen--Cahn & 1.50 & 5/5 & 1.0050 & 0.0042535 & 0.019222 \\
Allen--Cahn & 2.00 & 5/5 & 1.0047 & 0.0042535 & 0.032905 \\
\midrule
Cucker--Smale & 1.00 & 5/5 & 1.0674 & 0.00053947 & 0.005605 \\
Cucker--Smale & 1.25 & 5/5 & 1.2052 & 0.00053947 & 0.021064 \\
Cucker--Smale & 1.50 & 4/5 & 1.4951 & 47.2147 & 83.2411 \\
Cucker--Smale & 2.00 & 1/5 & 1.5540 & 352.3326 & 1782.6155 \\
\bottomrule
\end{tabular}
\end{table}

Table~\ref{tab:ood_long_horizon} shows that Allen--Cahn satisfies
the success criterion at all tested scales.
Cucker--Smale succeeds in all five cases through $s=1.25$, but
only four of five cases at $s=1.5$ and one of five
at $s=2$.
The larger final state norms and HJB residuals include
the failed cases and reveal the deterioration under
larger distribution shifts.
These results characterize empirical reliability
within the tested range.

\subsection{Sensitivity to Data Size and Dynamics}
\label{app:data_size_dynamics}

Both studies use the five $s=1$ initial conditions,
$T=50$, control-update interval of $0.01$, and success
criterion from Section~\ref{app:ood}.

\noindent\textbf{Data-size ablation.}
We compare the original full-data models with models
trained on half of the state--control pairs.
Subsets are sampled separately from the original
training and held-out partitions, preserving the
80/20 split for Cucker--Smale and the 90/10 split for Allen--Cahn. The half-data models retain the main-text training
settings in Section~\ref{sec:experiment}. With epoch counts and batch sizes fixed, halving the data
reduces the training workload, measured by optimizer updates,
from $500{,}000$ to $252{,}000$ for Cucker--Smale and from
$9{,}000$ to $6{,}000$ for Allen--Cahn.

\begin{table}[tbp]
\centering
\setlength{\tabcolsep}{5pt}
\caption{Data-size ablation.
Entries report mean cost ratio and mean
$E_{\rm AUC}$, defined in ~\eqref{eq:eauc},
over five common test initial conditions.
All five cases satisfy the success criterion
for each setting.}
\label{tab:data_size}
\begin{tabular}{lcrr}
\toprule
System & Train/held-out & Mean $C$ & Mean $E_{\rm AUC}$ \\
\midrule
Cucker--Smale & $16000/4000$ & 1.1210 & 0.0639 \\
Cucker--Smale & $8000/2000$  & 1.3742 & 0.1132 \\
\midrule
Allen--Cahn & $9000/1000$ & 1.0022 & 0.0324 \\
Allen--Cahn & $4500/500$  & 1.0027 & 0.0230 \\
\bottomrule
\end{tabular}
\end{table}

Table~\ref{tab:data_size} shows problem-dependent
sensitivity to data quantity.
Halving the data increases both cost and value--cost
discrepancy for Cucker--Smale, whereas Allen--Cahn retains near-reference
costs with a lower discrepancy in this run.
The reduced budget affects performance without
preventing successful regulation in these
fixed-seed evaluations.

\noindent\textbf{Dynamics perturbation.}
We freeze the original $\gamma=1$ models and their
input normalization, without retraining or modifying
their feedback laws.
For Cucker--Smale, we multiply the communication kernel
$P(y_i,y_j)$ used to construct the interaction matrix
in \eqref{eq:discrete_system_test1} by $0.9$, $1.0$, or $1.1$.
For Allen--Cahn, we replace the diffusion coefficient $0.1$
multiplying $A_0$ in ~\eqref{eq:discrete_system_test2}
by $0.09$, $0.10$, or $0.11$. All other dynamics and cost parameters remain unchanged.

The SDRE reference is recomputed using the actual
perturbed dynamics, and the HJB residual is evaluated
with those dynamics. All cases in Table~\ref{tab:dynamics_perturbation}
complete the full horizon and satisfy the regulation
criterion.
Cucker--Smale exhibits little change in median normalized cost.
For Allen--Cahn, reduced diffusion increases both normalized
cost and residual, although regulation remains successful.
These results support empirical tolerance to the
tested $\pm10\%$ single-parameter perturbations
within the specified initial-condition and
time-horizon settings.

\begin{table}[tbp]
\centering
\setlength{\tabcolsep}{4pt}
\caption{Frozen-policy evaluation under dynamics
perturbations.
The parameter is the communication-kernel multiplier
for Cucker--Smale and the diffusion coefficient for Allen--Cahn.
Costs are normalized by the corresponding
perturbed-system SDRE costs.
The cost column reports the median among successful
cases; terminal norms and absolute HJB residuals are
averaged over all five cases.}
\label{tab:dynamics_perturbation}
\begin{tabular}{lccrrr}
\toprule
System & Parameter & Success & Median $C$
& Mean $\|x(T)\|_2$ & Mean abs. HJB \\
\midrule
Cucker--Smale & 0.9 & 5/5 & 1.0693
& $5.3947\times10^{-4}$ & 0.005802 \\
Cucker--Smale & 1.0 & 5/5 & 1.0674
& $5.3947\times10^{-4}$ & 0.005605 \\
Cucker--Smale & 1.1 & 5/5 & 1.0680
& $5.3947\times10^{-4}$ & 0.005504 \\
\midrule
Allen--Cahn & 0.09 & 5/5 & 1.0247
& $6.1133\times10^{-3}$ & 0.010137 \\
Allen--Cahn & 0.10 & 5/5 & 1.0022
& $4.2535\times10^{-3}$ & 0.004296 \\
Allen--Cahn & 0.11 & 5/5 & 0.9984
& $3.2884\times10^{-3}$ & 0.007411 \\
\bottomrule
\end{tabular}
\end{table}


\end{document}